\documentclass[12pt]{article}
\usepackage{latexsym,amssymb,amsmath}
\usepackage{stmaryrd}
\begin{document}
\baselineskip=20pt

\newcommand{\la}{\langle}
\newcommand{\ra}{\rangle}
\newcommand{\psp}{\vspace{0.4cm}}
\newcommand{\pse}{\vspace{0.2cm}}
\newcommand{\ptl}{\partial}
\newcommand{\dlt}{\delta}
\newcommand{\sgm}{\sigma}
\newcommand{\al}{\alpha}
\newcommand{\be}{\beta}
\newcommand{\G}{\Gamma}
\newcommand{\gm}{\gamma}
\newcommand{\vs}{\varsigma}
\newcommand{\Lmd}{\Lambda}
\newcommand{\lmd}{\lambda}
\newcommand{\td}{\tilde}
\newcommand{\vf}{\varphi}
\newcommand{\yt}{Y^{\nu}}
\newcommand{\wt}{\mbox{wt}\:}
\newcommand{\rd}{\mbox{Res}}
\newcommand{\ad}{\mbox{ad}}
\newcommand{\stl}{\stackrel}
\newcommand{\ol}{\overline}
\newcommand{\ul}{\underline}
\newcommand{\es}{\epsilon}
\newcommand{\dmd}{\diamond}
\newcommand{\clt}{\clubsuit}
\newcommand{\vt}{\vartheta}
\newcommand{\ves}{\varepsilon}
\newcommand{\dg}{\dagger}
\newcommand{\tr}{\mbox{Tr}}
\newcommand{\ga}{{\cal G}({\cal A})}
\newcommand{\hga}{\hat{\cal G}({\cal A})}
\newcommand{\Edo}{\mbox{End}\:}
\newcommand{\for}{\mbox{for}}
\newcommand{\kn}{\mbox{ker}}
\newcommand{\Dlt}{\Delta}
\newcommand{\rad}{\mbox{Rad}}
\newcommand{\rta}{\rightarrow}
\newcommand{\mbb}{\mathbb}
\newcommand{\lra}{\Longrightarrow}
\newcommand{\X}{{\cal X}}
\newcommand{\Y}{{\cal Y}}
\newcommand{\Z}{{\cal Z}}
\newcommand{\U}{{\cal U}}
\newcommand{\V}{{\cal V}}
\newcommand{\W}{{\cal W}}
\newcommand{\sta}{\theta}
\setlength{\unitlength}{3pt}

\begin{center}{\Large \bf Representations of Lie Algebras }\end{center}
\begin{center}{\Large \bf and Coding Theory}\footnote
{2000 Mathematical Subject Classification. Primary 17B10, 94B60;
Secondary 17B25.}
\end{center}
\vspace{0.2cm}

\begin{center}{\large Xiaoping Xu}\end{center}
\begin{center}{Hua Loo-Keng Mathematical Laboratory}\end{center}
\begin{center}{Institute of Mathematics, Academy of Mathematics \& System Sciences}\end{center}
\begin{center}{Chinese Academy of Sciences, Beijing 100190, P.R. China}
\footnote{Research supported
 by China NSF 10871193}\end{center}

\begin{abstract}{ Linear codes with large minimal distances are important error correcting
codes in information theory.
Orthogonal codes have more applications in the other fields of mathematics. In this paper, we study the binary and
ternary orthogonal codes generated by the weight matrices on finite-dimensional  modules of
  simple Lie algebras. The Weyl groups of the Lie algebras act on these codes isometrically. It turns out
 that certain weight matrices of $sl(n,\mbb{C})$ and $o(2n,\mbb{C})$ generate doubly-even  binary orthogonal codes
 and ternary orthogonal codes with large minimal distances. Moreover, we prove that the weight matrices of $F_4$, $E_6$, $E_7$
 and $E_8$ on their minimal irreducible modules and adjoint modules all generate ternary orthogonal codes
 with large minimal distances. In determining the minimal distances, we have used the Weyl groups and branch rules of
 the irreducible representations of the related simple Lie algebras.}\end{abstract}

\section{Introduction}

Let $m$ be a positive integer and denote
$\mbb{Z}_m=\mbb{Z}/m\mbb{Z}$. A {\it code ${\cal C}$  of length $n$}
is a subset of $(\mbb{Z}_m)^n$ for some $m$, where the ring
structure of $\mbb{Z}_m$ may not be used. The elements of ${\cal C}$
are called {\it codewords}.  The {\it (Hamming) distance} between
two codewords is the number of different coordinates. The {\it
minimal distance} of a code is the minimal number among the
distances of all its pairs of codewords in the code. A code with
minimal distance $d$ can be used to correct
$\llbracket(d-1)/2\rrbracket$ errors in signal transmissions. We
refer [6], [15], [23] for more details. Examples of the well-known
infinite families of codes are cyclic codes, quadratic residue
codes, Goppa codes, algebraic geometry codes, arithmetic codes,
Hadamard codes and Pless double-circulant codes, etc. The names of
these families also indicate the methods of constructing codes. In
this paper, we introduce a new infinite family of codes arising from
finite-dimensional representations of simple Lie algebras, which we
may call {\it Lie theoretic codes}. One of the important features of
these codes is that the corresponding Weyl group acts on them
isometrically (may not be faithful).

A {\it linear code} ${\cal C}$ over the ring $\mbb{Z}_m$ is a
$\mbb{Z}_m$-submodule of $(\mbb{Z}_m)^n$. The {\it (Hamming) weight}
of a codeword in a linear code ${\cal C}$ is the number of its
nonzero coordinates. In this case, the minimal distance of ${\cal
C}$ is exactly the minimal weight of the nonzero codewords in ${\cal
C}$. The inner product in $(\mbb{Z}_m)^n$ is defined by
$$(a_1,...,a_n)\cdot(b_1,...,b_n)=\sum_{i=1}^na_ib_i.\eqno(1.1)$$
Moreover, ${\cal C}$ is called {\it orthogonal} if
$${\cal C}\subseteq\{\vec a\in (\mbb{Z}_m)^n\mid\vec\al\cdot\vec b=0\;\for\;\vec b\in{\cal C}\}.
\eqno(1.2)$$ When the equality holds, we call ${\cal C}$ a {\it
self-dual} code. Orthogonal linear codes (especially, self-dual
codes) have important applications to the other mathematical fields
such as sphere packings, integral linear lattices, finite group
theory, etc. We refer References [2]-[6], [9]-[14], [17]-[21] and
[24] for more details. A code is called {\it binary} if $m=2$ and
{\it ternary} when $m=3$. A binary linear code is called  {\it even
(doubly-even)} if the weights of all its codewords are divisible by
2 (by 4).

Let ${\cal G}$ be a finite-dimensional simple Lie algebras over $\mbb{C}$, the field of complex numbers.
Take a Cartan subalgebra $H$ and simple positive roots $\{\al_1,\al_2,...,\al_n\}$. Moreover, we denote by
$\{h_1,h_2,...,h_n\}$ the elements of $H$ such that the matrix
$$(\al_i(h_j))_{n\times n}\;\;\mbox{is the Cartan matrix of}\;\;{\cal G}\eqno(1.3)$$
(e.g., cf. [7]). For a finite-dimensional ${\cal G}$-module $V$,
it is well known that $V$ has a weight-subspace decomposition:
$$V=\bigoplus_{\mu\in H^\ast}V_\mu,\qquad V_\mu=\{v\in V\mid h(v)=\mu(h)v\;\for\;h\in H\}.\eqno(1.4)$$
Take a maximal linearly independent set $\{u_1,u_2,...,u_k\}$ of weight vectors with nonzero weights in $V$ such
 that the order is compatible with the partial order of weights (e.g., cf. [H]). Write
 $$h_i(u_j)=c_{i,j}u_j,\qquad C(V)=(c_{i,j})_{n\times k}.\eqno(1.5)$$
 By the representation theory of simple Lie algebras, all $c_{i,j}$ are integers. We call $C(V)$ the {\it weight matrix
 of ${\cal G}$ on} $V$. Identify  integers with their images in
 $\mbb{Z}_m$ when the context is clear. Denote by ${\cal C}_m(V)$ the linear code over $\mbb{Z}_m$ generated by
$C(V)$.

In this paper, we prove that ${\cal C}_2(V)$ and ${\cal C}_3(V)$
for certain finite-dimensional irreducible modules of special
linear Lie algebras are doubly-even binary orthogonal codes with
large minimal distances and ternary orthogonal codes with large
minimal distances, respectively. Moreover, ${\cal C}_3(V)$ for
certain finite-dimensional  modules of orthogonal Lie algebras are
also ternary orthogonal codes with large minimal distances.
Furthermore, we prove that the codes ${\cal C}_3(V)$ of the
exceptional simple Lie algebras $F_4$, $E_6$, $E_7$ and $E_8$ on
their minimal irreducible modules and adjoint modules are all
ternary orthogonal codes
 with large minimal distances. This coding theoretic phenomena was observed when we investigated the polynomial representations of these
 algebras in [26]-[28]. It is also well known that determining the minimal distance of a  linear
 code is in general very difficult. We have used the Weyl groups and branch rules of
 irreducible representations of the related simple Lie algebras in determining the minimal distances of the concerned
 codes. Note also that our code ${\cal C}_m(V)$ carries the important information of the simple root vectors acting on the weight vectors
 $u_i$ via the weight matrix $C(V)$ (e.g., c.f. [7]). Below we give more technical details.

 Suppose that the weight of $u_i$ is $\mu_i$. Set
 $${\cal H}_m=\sum_{i=1}^n\mbb{Z}_mh_i.\eqno(1.6)$$
We define a map $\Im: {\cal H}_m\rta (\mbb{Z}_m)^k$ by
$$\Im(\sum_{i=1}^nl_ih_i)=(\sum_{i=1}^nl_i\mu_1(h_i),\sum_{i=1}^nl_i\mu_2(h_i),...,\sum_{i=1}^nl_i\mu_k(h_i)).\eqno(1.7)$$
Then
$${\cal C}_m(V)=\Im({\cal H}_m).\eqno(1.8)$$
Denote by ${\cal W}({\cal G})$ the Weyl group of the simple Lie
algebra ${\cal G}$. For any $\sgm\in {\cal W}({\cal G})$, there
exists a linear automorphism $\hat\sgm$ of $V$ such that
$$\hat\sgm (V_{\mu})=V_{\sgm(\mu)},\qquad\sgm(\mu)(\sgm(h))=\mu(h)\qquad\for\;\;h\in H\eqno(1.9)$$
(e.g., cf. [7]). Moreover, we define an action of ${\cal W}({\cal
G})$ on ${\cal H}_m$ by
$$\sgm(\sum_{i=1}^nl_ih_i)=\sum_{i=1}^nl_i\sgm(h_i)\qquad\for\;\;\sgm\in {\cal W}({\cal G}).\eqno(1.10)$$
According to (1.9),
$$\wt\Im(\sgm(h))=\wt\Im(h)\qquad\for\;\;\sgm \in {\cal W}({\cal G}),\;h\in {\cal H}_m.\eqno(1.11)$$
So the number of the distinct weights of codewords in ${\cal
C}_m(V)$ is less than or equal to the number of ${\cal W}({\cal
G})$-orbits in ${\cal H}_m$. Expression (1.11) will be used later in
determining minimal distances.

Let $\Lmd(V)$ be the set of nonzero weights of $V$. The module $V$
is called {\it self-dual} if $\Lmd(V)=-\Lmd(V)$. In this paper, we
are only interested in the binary and ternary codes. We call ${\cal
C}_2(V)$ the {\it binary weight code of} ${\cal G}$ {\it on} $V$. If
$V$ is self-dual, then the weight matrix $C(V)=(-A,A)$ and ${\cal
C}_3(V)$ is orthogonal if and only if $A$ generates a ternary
orthogonal code (e.g., cf. [15]). For this reason, we call the
ternary code generated by $A$  the {\it ternary weight code of}
${\cal G}$ {\it on} $V$ if $V$ is self-dual. When $V$ is not
self-dual, then ${\cal C}_3(V)$ is the {\it ternary weight code of}
${\cal G}$ {\it on} $V$.

Denote by $V_X(\lmd)$ the finite-dimensional irreducible module of a
simple Lie algebra of type $X$ with the highest weight $\lmd$. Let
$p$ be a prime number. Then $\mbb{Z}_p$ is a finite field, which is
traditionally denoted as $\mbb{F}_p$. A linear code ${\cal C}$ of
length $n$ over $\mbb{F}_p$ is a linear subspace of $\mbb{F}_p^n$
over $\mbb{F}_p$. If $\dim {\cal C}=k$, we say that ${\cal C}$ is of
{\it type} $[n,k]$. When $d$ is the minimal distance of ${\cal C}$,
 we call ${\cal C}$ an {\it $[n,k,d]$-code.} Take the labels of
simple roots from [7]. Denote by $\lmd_i$ the $i$th fundamental
weight of the related simple Lie algebra. We summarize the main
results in this paper as the following three theorems.

The special linear Lie algebra $sl(n,\mbb{C})$ consists of all $n\times n$ matrices with zero trace, which is a simple Lie
algebra of type $A_{n-1}$.\psp

{\bf Theorem 1}. {\it (1) The binary weight code ${\cal
C}_2(V_{A_{2m-1}}(\lmd_2))$ of $sl(2m,\mbb{C})$ is a doubly-even
orthogonal $[m(2m-1),2(m-1),4(m-1)]$-code if $m\geq 2$.

(2) The binary weight code ${\cal C}_2(V_{A_{n-1}}(\lmd_3))$ of
$sl(n,\mbb{C})$ is a doubly-even  orthogonal $[{n\choose
3},n-1,(n-2)(n-3)]$-code if $n>9$ and $n\equiv
2,3\;(\mbox{mod}\;4)$.

(3)  The ternary weight code of $sl(3m+2,\mbb{C})$ on
$V_{A_{3m+1}}(\lmd_2)$  is an orthogonal  $[{3m+2\choose 2},\\
3m+1,6m]$-code if $m>0$.

(4)  The ternary weight code of $sl(3m,\mbb{C})$ on
$V_{A_{3m-1}}(\lmd_3)$ is is an orthogonal $[{3m\choose 3},3m-2,
3(m-1)(3m-2)]$-code. Moreover, the ternary weight code of
$sl(3m+2,\mbb{C})$ on $V_{A_{3m+1}}(\lmd_3)$ is  an orthogonal
$[{3m+2\choose 3},3m+1,3m(3m+1)/2]$-code.

(5) The ternary weight code of $sl(3m,\mbb{C})$ on the adjoint
module $sl(3m,\mbb{C})$
 is an orthogonal $[{3m\choose 2}, 3m-2,3(m-1)]$-code if $m>1$.}
\psp

The Lie algebra $o(2n,\mbb{C})$ consists of all $2n\times 2n$ skew-symmetric matrices, which is a simple
Lie algebra of type $D_n$.\psp

{\bf Theorem 2}. {\it
(1) The ternary weight code of $o(6m+2,\mbb{C})$ on $V_{D_{3m+1}}(\lmd_2)$ is
is an orthogonal  $[2m(3m+1),3m+1,6m]$-code if $m>0$.

(2)  The ternary weight code of $o(2m,\mbb{C})$ on $V_{D_m}(\lmd_3)$ is
is an orthogonal $[m(m-1)(2m-1)/3,m,(m-1)(2m-3)]$-code if
 $m\not\equiv -1\;(\mbox{mod}\;3)$ and $m>3$.

(3) The ternary code ${\cal C}_3(V_{D_m}(\lmd_m))$ of
$o(2m,\mbb{C})$ is of type $[2^{m-1},m,2^{m-2}]$ if $6\neq m>3$
and of type $[32,6,12]$ when $m=6$, where the representation of
$o(2m,\mbb{C})$ on ${\cal C}_3(V_{D_m}(\lmd_m))$ is the spin
representation.

(4)  The ternary weight code of $o(12m+4,\mbb{C})$ on
$o(12m+4,\mbb{C})+V_{D_{6m+2}}(\lmd_{6m+2})$ is an orthogonal
$[(6m+2)(6m+1)+2^{6m},6m+2,24m+1+2^{6m-1}]$-code for $m>0$.} \psp

There are five exceptional finite-dimensional simple Lie algebras,
labeled as $G_2,\;F_4,\;E_6,$ $E_7$ and $E_8$. They have broad
applications. We find the following common coding theoretic feature
of the simple Lie algebras of types $F_4,\;E_6,\;E_7$ and $E_8$.

\psp

{\bf Theorem 3}. {\it (1) The ternary weight code of $F_4$ on its minimal module is an orthogonal [12,4,6]-code.

(2) The ternary weight code of $F_4$ on its adjoint module is an orthogonal [24,4,15]-code.

(3) The ternary weight code of $E_6$ on its minimal module is an orthogonal [27,6,12]-code.

(4) The ternary weight code of $E_6$ on its adjoint module is an orthogonal [36,5,21]-code.

(5) The ternary weight code of $E_7$ on its minimal module is an orthogonal [28,7,12]-code.

(6) The ternary weight code of $E_7$ on its adjoint module is an orthogonal [63,7,27]-code.

(7) The ternary weight code of $E_8$ on its minimal (adjoint)
module is an orthogonal [120,8,57]-code.} \psp

Section 2 is devoted to the study of the binary and ternary weight
codes of  $sl(n,\mbb{C})$. In Section 3, we prove Theorem 2.
Section 4 is about the ternary weight codes of $F_4$ on its
minimal module and adjoint module. In Section 5, we investigate
the ternary weight codes of $E_6$ on its minimal module and
adjoint module. We deal with the ternary weight codes of $E_7$ and
$E_8$ on their minimal module and adjoint module in Section 6.

\section{Codes Related to Representations of $sl(n,\mbb{C})$}

In this section, we study the binary and ternary codes related to
representations of $sl(n,\mbb{C})$, where $n>1$ is an integer.

Throughout this paper, we always take the following notion:
$$\ol{i,j}=\{i,i+1,i+2,...,j\}\eqno(2.1)$$
for any integers $i\leq j$. We denote
$$\ves_i=(0,...,\stl{i}{1},0,...,0)\in\mbb{R}^n.\eqno(2.2)$$
So
$$\mbb{R}^n=\sum_{i=1}^n\mbb{R}\ves_i.\eqno(2.3)$$
Then inner product ``$(\cdot,\cdot)$" is Euclidian, that is,
$$(\sum_{i=1}^nk_i\ves_i,\sum_{j=1}^nl_j\ves_j)=\sum_{i=1}^nk_il_i.\eqno(2.4)$$

Denote by $E_{i,j}$ the square matrix with 1 as its $(i,j)$-entry
and 0 as the others.  The special linear Lie algebra
$$sl(n,\mbb{C})=\sum_{n\leq i<j\leq
n}(\mbb{C}E_{i,j}+\mbb{C}E_{j,i})+\sum_{r=1}^{n-1}\mbb{C}h_r,\qquad
h_r=E_{r,r}-E_{r+1,r+1}.\eqno(2.5)$$ The subspace
$H_{A_{n-1}}=\sum_{i=1}^{n-1}\mbb{C}h_i$ forms a Cartan subalgebra
of $sl(n,\mbb{C})$. The root system
$$\Phi_{A_{n-1}}=\{\ves_i-\ves_j\mid i,j\in\ol{1,n},\;i\neq
j\}.\eqno(2.6)$$ Take the simple positive roots
$$\al_i=\ves_i-\ves_{i+1}\qquad\for\;\;i\in\ol{1,n-1}.\eqno(2.7)$$
The corresponding Dynkin diagram is\psp

\begin{picture}(70,8)\put(2,0){$A_{n-1}$:}
\put(21,0){\circle{2}}\put(22,0){\line(1,0){12}}\put(21,-5){1}
\put(35,0){\circle{2}}\put(35,-5){2}\put(43,0){...}\put(53,0){\circle{2}}
\put(53,-5){n-2}
\put(54,0){\line(1,0){12}}\put(67,0){\circle{2}}\put(67,-5){n-1}
\end{picture}
\vspace{1cm}

\noindent The Weyl group ${\cal W}_{A_{n-1}}$ of $sl(n,\mbb{C})$
is exactly the full permutation group $S_n$ on $\ol{1,n}$, which
acts on $H_{A_{n-1}}$ and $\mbb{R}^n$ by permuting sub-indices of
$E_{i,i}$ and $\ves_i$, respectively.

Let ${\cal A}$ be the associative algebra generated by
$\{\sta_1,\sta_2,...,\sta_n\}$ with the defining relations:
$$\sta_i\sta_j=-\sta_j\sta_i\qquad\for\;\;i,j\in\ol{1,n}.\eqno(2.8)$$
The generators $\sta_i$ are called {\it spin variables}. The
representation of the Lie algebra $sl(n,\mbb{C})$ on ${\cal A}$ is
given by
$$E_{i,j}=\sta_i\ptl_{\sta_j}\qquad\for\;\;i,j\in\ol{1,n}.\eqno(2.9)$$
Set
$${\cal A}_r=\sum_{1\leq i_1<i_2<\cdots< i_r\leq
n}\mbb{C}\sta_{i_1}\sta_{i_2}\cdots\sta_{i_r}\qquad\for\;\;r\in\ol{1,n}.
\eqno(2.10)$$ Then ${\cal A}_r$ forms an irreducible
$sl(n,\mbb{C})$-submodule of highest weight $\lmd_r$ for
$r\in\ol{1,n-1}$, that is, ${\cal A}_r\cong V_{A_{n-1}}(\lmd_r)$.
The Weyl group ${\cal W}_{A_{n-1}}$ acts on ${\cal A}$ by permuting
sub-indices of $\sta_i$.

Two $k_1\times k_2$ matrices $A_1$ and $A_2$ with entries in
$\mbb{Z}_m$ are called {\it equivalent} in the sense of coding
theory if there exist an invertible $k_1\times k_1$ matrix $K_1$ and
an invertible $k_2\times k_2$ monomial matrix $K_2$ such that
$A_1=K_1A_2K_2.$ Equivalent matrices generate isomorphic codes.
 Take any order of the basis
$$\{x_{r,1},x_{r,2},...,x_{r,{n\choose r}}\}=\{
\sta_{i_1}\sta_{i_2}\cdots\sta_{i_r}\mid 1\leq i_1<i_2<\cdots<
i_r\leq n\}.\eqno(2.11)$$ Then we have
$$h_i(x_{r,j})=a_{i,j}(r)x_{r,j},\qquad
a_{i,j}(r)\in\mbb{Z}.\eqno(2.12)$$ Modulo equivalence, the weight
matrix
$$C({\cal A}_r)=[a_{i,j}(r)]_{(n-1)\times {n\choose
r}}.\eqno(2.13)$$

{\bf Theorem 2.1}. {\it When $n=2m\geq 4$ is even, ${\cal C}_2({\cal
A}_2)$ is a doubly-even binary orthogonal
$[m(2m-1),2(m-1),4(m-1)]$-code}.

{\it Proof}. Denote by $\xi_i$ the $i$th row $C_2({\cal A}_2)$. Then
$$\wt\xi_i=2(n-2)\qquad\for\;\;i\in\ol{1,n-1}.\eqno(2.14)$$
Moreover,
$$\sum_{i=0}^{m-1}\xi_{2i+1}=0\qquad\mbox{in}\;\;{\cal C}_2({\cal
A}_2). \eqno(2.15)$$ Furthermore,
$$\xi_i\cdot\xi_j=4\equiv 0\qquad \mbox{if}\;\;i+1<j\eqno(2.16)$$
and
$$\xi_i\cdot\xi_{i+1}=2(m-1)\equiv 0.\eqno(2.17)$$

Write
$$E_{i,i}(x_{r,j})=b_{i,j}(r)x_{r,j},\qquad
B_r=[b_{i,j}(r)]_{n\times {n\choose r}}.\eqno(2.18)$$ Denote by
$\zeta_i$ the $i$th row of $B_2$. By symmetry (cf. (1.9)-(1.11)),
any nonzero codeword in ${\cal C}_2({\cal A}_2)$ has the same
weight as the codeword
$$u=\sum_{s=1}^{2t}\zeta_s\in\mbb{F}_2^{n(n-1)/2}\qquad\mbox{for
some}\;\;t\in\ol{1,m-1}.\eqno(2.19)$$ We calculate
$$\wt u=4t(m-t)=-4t^2+4mt.\eqno(2.20)$$ Since the function $-4t^2+t(4m-1)$ attains
maximal at $t=m/2$, $\wt u$ is minimal at $t=1$ or $m-1$. Note
$$\wt u=4(m-1)\qquad\mbox{if}\;\;t=1\;\mbox{or}\;m-1.\eqno(2.21)$$
Thus the code ${\cal C}_2({\cal A}_2)$ has the minimal distance
$4(m-1).\qquad\Box$\psp

When $m=2$, ${\cal C}_2({\cal A}_2)$ is a doubly-even binary
orthogonal $[6,2,4]$-code. If $m=3$,
 ${\cal C}_2({\cal A}_2)$ becomes a doubly-even binary orthogonal $[15,4,8]$-code.
These two code are optimal linear
 codes (e.g., cf. [1]).
In  the case of $m=4$,  ${\cal C}_2({\cal A}_2)$ is a doubly-even
binary orthogonal $[28,6,12]$-code.\psp

{\bf Theorem 2.2}. {\it The code ${\cal C}_2({\cal A}_3)$ is a
doubly-even binary orthogonal $[{n\choose 3}, n-1,(n-2)(n-3)]$-code
if $n>9$ and $n\equiv 2,3\;(\mbox{mod}\;4)$.}

{\it Proof}. Denote by $\xi_i$ the $i$th row the weight matrix
$C({\cal A}_3)$. Then
$$\wt\xi_i=(n-2)(n-3)\qquad\for\;\;i\in\ol{1,n-1}.\eqno(2.22)$$
Moreover,
$$\xi_i\cdot\xi_j=4(n-4)\qquad \mbox{if}\;\;i+1<j\eqno(2.23)$$
and
$$\xi_i\cdot\xi_{i+1}=n-3+{n-3\choose 2}=\frac{(n-2)(n-3)}{2}.\eqno(2.24)$$
So ${\cal C}_2({\cal A}_3)$ is a doubly-even binary orthogonal code
under the assumption.

Denote by $\zeta_i$ the $i$th row of $B_3$ (cf. (2.18)). By symmetry
(cf. (1.9)-(1.11)), any nonzero codeword in ${\cal C}_2({\cal A}_3)$
has the same weight as the codeword
$$u(t)=\sum_{s=1}^{2t}\zeta_s\in\mbb{F}_2^n\qquad\mbox{for
some}\;\;t\in\ol{1,\llbracket n/2\rrbracket}.\eqno(2.25)$$ We
calculate
$$f(t)=3\wt u(t)=3{2t\choose 3}+6t{n-2t\choose
2}=t[16t^2-12nt+3n(n-1)+2].\eqno(2.26)$$ Moreover,
$$f'(t)=48t^2-24nt+3n(n-1)+2=48\left(t-\frac{n}{4}\right)^2-3n+2.\eqno(2.27)$$
Thus
$$f'(t_0)=0\lra
t_0=\frac{n}{4}\pm\frac{1}{4}\sqrt{n-\frac{2}{3}}.\eqno(2.28)$$
Since $f'(0)=3n(n-1)+2>0$, $f(t)$ attains local maximum at
$$t=\frac{n}{4}-\frac{1}{4}\sqrt{n-\frac{2}{3}}\eqno(2.29)$$
and local minimum at
$$t=\frac{n}{4}+\frac{1}{4}\sqrt{n-\frac{2}{3}}.\eqno(2.30)$$
According to (2.22) and (2.26), $f(1)=3(n-2)(n-3)$. Furthermore,
\begin{eqnarray*}& &f\left(\frac{n}{4}+\frac{1}{4}\sqrt{n-\frac{2}{3}}\right)
\\ &=&\left(\frac{n}{4}+\frac{1}{4}\sqrt{n-\frac{2}{3}}\right)
\left[16\left(\frac{n}{4}+\frac{1}{4}\sqrt{n-\frac{2}{3}}\right)^2
-12n\left(\frac{n}{4}+\frac{1}{4}\sqrt{n-\frac{2}{3}}\right)
+3n(n-1)+2\right]\\&=&\frac{1}{4}\left(n+\sqrt{n-\frac{2}{3}}\right)
\left[\left(n+\sqrt{n-\frac{2}{3}}\right)^2-3n\left(n+
\sqrt{n-\frac{2}{3}}\right)
+3n(n-1)+2\right]\\&=&\frac{1}{4}\left(n+\sqrt{n-\frac{2}{3}}\right)
\left[n\left(n-\sqrt{n-\frac{2}{3}}\right)-2n+\frac{4}{3}\right]\\
&=&\frac{1}{4}\left[n^3-3n^2+2n+\left(\frac{4}{3}-2n\right)
\sqrt{n-\frac{2}{3}}\right]\\
&>&\frac{1}{4}(n^3-5n^2+2n).\hspace{10.6cm}(2.31)\end{eqnarray*}
Thus
\begin{eqnarray*}& &
f\left(\frac{n}{4}+\frac{1}{4}\sqrt{n-\frac{2}{3}}\right)-f(1)\\
&>&\frac{1}{4}(n^3-5n^2+2n)-3(n-2)(n-3)=\frac{1}{4}(n^3-17n^2+62n-72)
\\ &>&\frac{n^2(n-17)}{4}.\hspace{11.7cm}(2.32)\end{eqnarray*}
If $n\geq 17$, we have
$$f\left(\frac{n}{4}+\frac{1}{4}\sqrt{n-\frac{2}{3}}\right)>f(1)
\eqno(2.33)$$ and \begin{eqnarray*}\hspace{2cm}& &f(n/2)-f(1)\\&=&
\frac{n}{2}[4n^2-6n^2+3n(n-1)+2]-3(n-2)(n-3)
\\ &=&\frac{n(n-1)(n-2)}{2}-3(n-2)(n-3)\\&=&\frac{(n-2)(n^2-7n+9)}{2}>0\;\;\mbox{if}\;\;n\geq 6.
\hspace{5.6cm}(2.34)\end{eqnarray*} Thus the minimal weight is
$f(1)/3=(n-2)(n-3)$ when $n\geq 17$.

When $n=10$, we calculate
\begin{center}{\bf Table 2.1}\end{center}
\begin{center}\begin{tabular}{|c|c|c|c|c|c|}\hline
$t$&1&2&3&4&5
\\\hline\wt u(t)&56&64&56&64&120\\\hline\end{tabular}\end{center}
If $n=11$, we find
\begin{center}{\bf Table 2.2}\end{center}
\begin{center}\begin{tabular}{|c|c|c|c|c|c|}\hline
$t$&1&2&3&4&5
\\\hline\wt u(t)&72&88&80&80&120\\\hline\end{tabular}\end{center}
When $n=14$, we obtain
\begin{center}{\bf Table 2.3}\end{center}
\begin{center}\begin{tabular}{|c|c|c|c|c|c|c|c|}\hline
$t$&1&2&3&4&5&6&7
\\\hline\wt u(t)&132&184&188&176&180&232&364\\\hline\end{tabular}\end{center}
If $n=15$, we get
\begin{center}{\bf Table 2.4}\end{center}
\begin{center}\begin{tabular}{|c|c|c|c|c|c|c|c|}\hline
$t$&1&2&3&4&5&6&7
\\\hline\wt u(t)&156&224&216&224&220&256&364\\\hline\end{tabular}\end{center}
This prove the conclusion in the theorem.$\qquad\Box$\psp

Note that when $n=6$, we find \begin{center}{\bf Table
2.5}\end{center}
\begin{center}\begin{tabular}{|c|c|c|c|}\hline
$t$&1&2&3
\\\hline\wt u(t)&12&8&20\\\hline\end{tabular}\end{center}
So ${\cal C}_2({\cal A}_3)$ is a doubly-even binary orthogonal
$[20,5,8]$-code. Moreover, if $n=7$, we find
\begin{center}{\bf Table 2.6}\end{center}
\begin{center}\begin{tabular}{|c|c|c|c|}\hline
$t$&1&2&3
\\\hline\wt u(t)&20&16&20\\\hline\end{tabular}\end{center}
 Hence ${\cal C}_2({\cal A}_3)$ a doubly-even binary orthogonal
$[35,6,16]$-code. In both cases, the above theorem fails and both
codes are the best even codes among the binary codes with the same
length and dimension (e.g., cf. [1]).

  According
to the above theorem, ${\cal C}_2({\cal A}_3)$ is a doubly-even
binary orthogonal $[120,9,56]$-code when $n=10$, $[165,10,72]$-code
if $n=11$, $[364,13,132]$-code when $n=14$ and
 $[455,14,156]$-code if $n=15$.

Next let us consider the ternary codes.  Again by symmetry, any
nonzero codeword in ${\cal C}_3({\cal A}_r)$ has the same weight as
the codeword
$$u(s,t)=\sum_{r=1}^s\zeta_r-\sum_{i=1}^t\zeta_{s+i}\in\mbb{F}_3^{{n\choose r}}\eqno(2.35)$$
for some nonnegative integers $s,t$, where $\zeta_\iota$ is the
$\iota$th row of the matrix $B_r$ in (2.18).
 Moreover,
$$\wt u(s,t)=\wt u(t,s).\eqno(2.36)$$
Furthermore, we have
$$\wt u(s,t)=(s+t)(n-s-t)+{s\choose 2}+{t\choose 2}\qquad\mbox{in}\;\;
{\cal C}_3({\cal A}_2)\eqno(2.37)$$ and
$$\wt u(s,t)=(s+t){n-s-t\choose 2}+(n-s){s\choose 2}+(n-t){t\choose 2}\qquad\mbox{in}\;\;
{\cal C}_3({\cal A}_3).\eqno(2.38)$$

For convenience, we denote \begin{eqnarray*}
\hspace{1cm}f(s,t)&=&2\wt u(s,t)=2(s+t)(n-s-t)+s(s-1)+t(t-1)\\
&=&(2n-1)(s+t)-s^2-t^2-4st\hspace{6.4cm}(2.39)\end{eqnarray*} in
${\cal C}_3({\cal A}_2)$ and
 \begin{eqnarray*}& &g(s,t)\\ &=&2\wt u(s,t)=
(s+t)(n-s-t)(n-s-t-1)+(n-s)s(s-1)+(n-t)t(t-1)
\\
&=&(s+t)^3-(2n-1)(s+t)^2+n(n-1)(s+t)-s^3-t^3
+(n+1)(s^2+t^2)-n(s+t)
\\&=&3st^2+3s^2t+(2-n)(s^2+t^2)-2(2n-1)st+n(n-2)(s+t)\hspace{2.9cm}(2.40)\end{eqnarray*} in
${\cal C}_3({\cal A}_3)$.

Note
$$f(3,0)=3(2n-1)-9=6(n-2),\;\;f(n,0)=n(2n-1)-n^2=n(n-1),\eqno(2.41)$$
$$f(1,1)=2(2n-1)-6=4(n-2),\;\;f(1,n-1)=(n-1)(n-2).\eqno(2.42)$$
Since geometrically $f(s,t)$ has only local minimum, it attains the
absolute minimum at boundary points. Thus
$$\min \{f(s,t)\mid s\equiv t\;(\mbox{mod}\;3)\}=4(n-2)\qquad\mbox{if}\;\;n\geq 5.\eqno(2.43)$$

Now
$$g_s(s,t)=3t^2+6st+2(2-n)s-2(2n-1)t+n(n-2),\eqno(2.44)$$
$$g_t(s,t)=3s^2+6st+2(2-n)t-2(2n-1)s+n(n-2). \eqno(2.45)$$ Suppose that
$g_s(s_0,t_0)=g_t(s_0,t_0)=0$ for $s_0,t_0\geq 0$, that is,
$$3t_0^2+6s_0t_0+2(2-n)s_0-2(2n-1)t_0+n(n-2)=0,\eqno(2.46)$$
$$3s_0^2+6s_0t_0+2(2-n)t_0-2(2n-1)s_0+n(n-2)=0.\eqno(2.47)$$
By $(2.46)-(2.47)$, we get
$$(t_0-s_0)(3t_0+3s_0-2(n+1))=0\lra t_0=s_0\;\;\mbox{or}\;\;
3t_0+3s_0=2(n+1).\eqno(2.48)$$

If $s_0=t_0$, then
 we find
$$9s_0^2-2(n-1)s_0+n(n-2)=0\sim
8s_0^2+(s_0-n+1)^2-1=0,\eqno(2.49)$$ which leads to a
contradiction because $n>1$. Thus $3t_0+3s_0=2(n-1)$. Denote
$s_1=3t_0$ and $t_1=3t_0$. Then $s_1+t_1=2(n+1)$ and (2.46)
becomes
$$t_1^2+2(2(n+1)-t_1)t_1+2(2-n)(2(n+1)-t_1)-2(2n-1)t_1+3n(n-2)=0,
\eqno(2.50)$$ equivalently,
$$t_1^2-2(n+1)t_1+(n-2)(n+4)=0\sim (t_1-n-1)^2-9=0
\lra t_1=n+4,\;n-2.\eqno(2.51)$$ Therefore,
$$s_0=\frac{n+4}{3},\;\;t_0=\frac{n-2}{3}\qquad\mbox{or}\qquad
t_0=\frac{n+4}{3},\;\;s_0=\frac{n-2}{3}.\eqno(2.52)$$ We calculate
$$g(s_0,t_0)=\frac{2(n-2)(n^2-n-3)}{9},\eqno(2.53)$$
$$g(1,0)=g(n-1,0)=(n-1)(n-2),\;\;
g(3,0)=3(n-2)(n-3),\;\;g(n,0)=0.\eqno(2.54)$$
$$g(1,1)=g(n-2,1)=2(n-2)(n-3),\;\;g(n-2,0)=2(n-2)^2.\eqno(2.55)$$
Moreover,
$$g(s_0,t_0)\geq g(1,0),\;g(1,1)\qquad\mbox{if}\;\;n\geq 6.\eqno(2.56)$$
When $n=5$, we calculate
$$g(1,0)=g(1,1)=g(2,1)=g(2,2)=g(3,1)=g(4,0)=g(4,1)=12,\eqno(2.57)$$
$$g(2,0)=g(3,0)=g(3,2)=18.\eqno(2.58)$$
In summary, we have:\psp

{\bf Theorem 2.3}. {\it Let $n\geq 5$. The the matrix $B_3$ (cf.
(2.18)) generates a ternary $\left[{n\choose 3},n-1,{n-1\choose
2}\right]$-code, which is equal to ${\cal C}_3({\cal A}_3)$ if
$n\not\equiv 0\;(\mbox{\it mod}\;3)$. If $n=3m+2$ for some positive
integer $m$, ${\cal C}_3({\cal A}_2)$ is a ternary orthogonal
$[{3m+2\choose 2},3m+1,6m]$-code and ${\cal C}_3({\cal A}_3)$ is a
ternary orthogonal $[{3m+2\choose 3},3m+1,3m(3m+1)/2]$-code.  The
code ${\cal C}_3({\cal A}_3)$ is a ternary orthogonal
$\left[{n\choose 3},n-2,(n-2)(n-3)\right]$-code when $n\equiv
0\;(\mbox{\it mod}\;3)$. }

{\it Proof}. The part of minimal distances has been proved by the
above arguments. We only need to prove orthogonality.

Suppose $n=3m+2$. In ${\cal C}_3({\cal A}_2)$, $\xi_r$ stands for
the $r$th row of  the weight matrix $C({\cal A}_2)$ and
$$\xi_i\cdot \xi_j=2-2=0\qquad\for\;\;1\leq i<j-1\leq
n-2,\eqno(2.59)$$
$$\xi_r\cdot
\xi_{r+1}=-(n-2)=-3m,\;\;\xi_s\cdot\xi_s=2(n-2)=6m\eqno(2.60)$$ for
$r\in\ol{1,n-2}$ and $s\in\ol{1,n-1}$. So ${\cal C}_3({\cal A}_2)$
is orthogonal. Now $\zeta_r$ stands for the $r$th row of $B_3$ (cf.
(2.18)). Observe
$$\sum_{i=1}^n\zeta_i=0\in \mbb{F}_3^{{n\choose 3}}\eqno(2.61)$$
by (2.9) and (2.10). Moreover,
$$\zeta_i\cdot\zeta_j=n-2=3m,\;\;\zeta_i\cdot\zeta_i={n-1\choose
2}=\frac{3m(3m+1)}{2},\qquad i\neq j.\eqno(2.62)$$ Thus $B_3$
generate a ternary  orthogonal code.

Assume that $n=3m$ for some nonnegative integer $m$. In ${\cal
C}_3({\cal A}_3)$, we also use $\xi_r$ for the $r$th row of the
weight code $C({\cal A}_3)$ and
$$\xi_i\cdot\xi_j=2(n-4)-2(n-4)=0\qquad\for\;\;1\leq i<j-1\leq
n-2,\eqno(2.63)$$
$$\xi_s\cdot\xi_s=2\xi_r\cdot\xi_{r+1}=(n-2)(n-3)=3(3m-2)(m-1)\equiv
0\eqno(2.64)$$ for $r\in\ol{1,n-2}$ and $s\in\ol{1,n-1}$. So ${\cal
C}_3({\cal A}_3)$ is orthogonal.$\qquad\Box$ \psp

According to the above theorem, ${\cal C}_3({\cal A}_2)$ is a
ternary orthogonal $[10,4,6]$-code when $n=5$ (which is optimal
(e.g., cf. [1])), $[28,7,12]$-code when $n=8$, and $[55,10,18]$-code
when $n=11$. Moreover, ${\cal C}_3({\cal A}_3)$ is a ternary
orthogonal $[10,4,6]$-code when $n=5$, $[15,4,12]$-code if $n=6$,
$[56,7,21]$-code when $n=8$, $[84,7,42]$-code if $n=9$,
$[165,10,45]$-code when $n=11$ and $[220,10,90]$-code when $n=12$.

Finally, we consider the adjoint representation of
$sl(n,\mbb{C})$. Note that $\{E_{i,j}\mid 1\leq i< j\leq n\}$ are
positive root vectors. Given an order
$$\{y_1,...,y_{{n\choose 2}}\}=\{E_{i,j}\mid 1\leq i< j\leq
n\},\eqno(2.65)$$ we write
$$[h_i,y_j]=k_{i,j}y_j,\qquad
[E_{r,r},y_j]=l_{r,j}y_j.\eqno(2.66)$$  Denote
$$K=(k_{i,j})_{(n-1)\times {n\choose r}},
\qquad L=(l_{i,j})_{n\times {n\choose r}}.\eqno(2.67)$$ Let ${\cal
K}$ be the ternary code generated by $K$ and let ${\cal L}$ be the
ternary code generated by $L$. Moreover, $\vec k_i$ stands for the
$i$th row of $K$ and $\vec l_r$ stands for the $r$th row of $L$. Set
$$u(s,t)=\sum_{i=1}^s\vec l_i-\sum_{j=1}^t\vec
l_{s+j}.\eqno(2.68)$$ For any nonzero codeword $v\in{\cal L}$, using
negative root vectors, we can prove
$$\wt (v,-v)=\wt (u(s,t),-u(s,t))\eqno(2.69)$$
for some $s$ and $t$ by symmetry (cf. (1.9)-(1.11)). Thus
$$\wt v=\wt u(s,t)=(s+t)(n-s-t)+st=\phi(s,t).\eqno(2.70)$$
Note
$$\phi(s,t)=n^2-\frac{1}{2}[(s-n)^2+(t-n)^2+(s-t)^2].\eqno(2.71)$$
So $\phi(s,t)$ has only local maximum. Thus it attains the
absolute minimum at the boundary points. We calculate
$$\phi(1,0)=\phi(n-1,0)=n-1,\;\;\phi(n-3,0)=3(n-3),\eqno(2.72)$$
$$\phi(1,1)=2n-3,\qquad\phi(n-2,1)=2(n-1).\eqno(2.73)$$
Since
$$\sum_{i=1}^n\vec l_i=0,\eqno(2.74)$$
$${\cal K}={\cal L}\qquad\mbox{if}\;\;n\neq
0\;(\mbox{mod}\;3).\eqno(2.75)$$
$$\vec k_i\cdot\vec k_j=2-2=0\qquad1\leq i<j-1\leq n,\eqno(2.76)$$
$$\vec k_r\cdot\vec k_{r+1}=6-n,\;\;\vec k_s\cdot\vec
k_s=2n-3.\eqno(2.77)$$ In summary, we have:\psp

{\bf Theorem 2.4}. {\it The code ${\cal L}$ is a ternary $[{n\choose
2},n-1,n-1]$-code if $n\geq 4$, which is also the ternary weight
code on the adjoint module $sl(n,\mbb{C})$ when $n\neq
0\;(\mbox{mod}\;3)$. If $n=3m$ for some integer $m>1$, then the
ternary weight code ${\cal K}$ on $sl(3m,\mbb{C})$ is an orthogonal
$[{3m\choose 2}, 3m-2, 3(m-1)]$-code.}

\section{Codes Related to Representations of $o(2m,\mbb{C})$}

In this section, we only study ternary codes related to certain
representations of $so(2m,\mbb{C})$, some of which will be used to
investigate the codes related to exceptional  simple Lie algebras.

 Let $n=2m$ be a positive even integer. Take the settings in (2.1)-(2.4) (with $n\rta m$).  The
orthogonal Lie algebra
\begin{eqnarray*}o(2m,\mbb{C})&=&\sum_{1\leq i<j\leq
m}[\mbb{C}(E_{i,j}-E_{m+j,m+i})+\mbb{C}(E_{j,i}-E_{m+i,m+j})
+\mbb{C}(E_{i,m+j}-E_{j,m+i})\\ & &+\mbb{C}(E_{m+i,j}-E_{m+j,i})]
+\sum_{r=1}^m\mbb{C}h_r,\hspace{6.5cm}(3.1)\end{eqnarray*} where
$$h_s=E_{s,s}-E_{s+1,s+1}-E_{m+s,m+s}-E_{m+s+1,m+s+1}\qquad\for\;\;s
\in\ol{1,m-1}\eqno(3.2)$$ and
$$h_m=E_{m-1,m-1}+E_{m,m}-E_{2m-1,2m-1}-E_{2m,2m}.\eqno(3.3)$$
Indeed, we take the Cartan subalgebra
$$ H_{D_m}=\sum_{i=1}^m\mbb{C}h_i\eqno(3.4)$$
of $o(2m,\mbb{C})$.  The root system
$$\Phi_{D_m}=\{\pm\ves_i\pm\ves_j\mid i,j\in\ol{1,m},\;i\neq
j\}\eqno(3.5)$$ and simple positive roots are:
$$\al_i=\ves_i-\ves_{i+1},\;\;\al_m=\ves_{m-1}+\ves_m,\qquad
i\in\ol{1,m-1}.\eqno(3.6)$$ The corresponding Dynkin diagram is

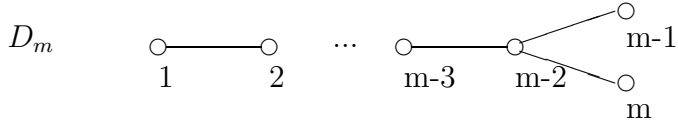
\begin{picture}(81,15)\put(2,0){$D_m$}\put(21,0){\circle{2}}\put(21,-5){1}\put(22,0){\line(1,0)
{12}} \put(35,0){\circle{2}}\put(35,-5){2}
\put(43,0){...}\put(52,0){\circle{2}}\put(52,-5){m-3}\put(53,0){\line(1,0){12}}
\put(66,0){\circle{2}}\put(66,-5){m-2}\put(66.6,0.5)
{\line(3,1){12}}\put(80,4.5){\circle{2}}\put(80,0){m-1}\put(66.6,-0.5)
{\line(3,-1){12}}\put(80,-4.5){\circle{2}}\put(80,-9){m}
\end{picture}
\vspace{1cm}

The Weyl group is $S_m\ltimes\mbb{Z}_2^{m-1}$, which acts $H_{D_m}$
and $\mbb{R}^m$ by permuting sub-indices of $\ves_i$ and
$E_{i,i}-E_{m+i,m+i}$, and changing sign on even number of their
coefficients.

Take the settings in (2.8)-(2.13) and (2.18).  Moreover, the
representation of
 $o(2m,\mbb{C})$ on ${\cal A}$ determined by (2.9).
 For any $\vec\iota=(\iota_1,...,\iota_m)$ with $\iota_i\in\{0,1\}$ and $\tau\in S_m$, we have an associative algebra
automorphism $\sgm_{\tau,\vec\iota}$ of ${\cal A}$ determined by
$$\sgm_{\tau,\vec\iota}(\sta_i)=
\sta_{m\dlt_{\iota_i,1}+\tau(i)},\;\;
\sgm_{\tau,\vec\iota}(\sta_{m+i})=\sta_{m\dlt_{\iota_i,0}+\tau(i)}
\qquad\for\;\;i\in\ol{1,m}.\eqno(3.6)$$ Moreover, we define a linear
map $\sgm_{\tau,\vec\iota}$ on ${\cal H}$ by
$$\sgm_{\tau,\vec\iota}(E_{i,i}-E_{m+i,m+i})=(-1)^{\iota_i}(E_{\tau(i),\tau(i)}-E_{m+\tau(i),m+\tau(i)})\qquad\for
\;\;i\in\ol{1,m}.\eqno(3.7)$$ Then
$$\sgm_{\tau,\vec\iota}(h(w))=\sgm_{\tau,\vec\iota}(h)[\sgm_{\tau,\vec\iota}(w)]\qquad\for\;\;h\in{\cal H},\;
w\in{\cal A}.\eqno(3.8)$$

Note that all ${\cal A}_r\cong V_{D_m}(\lmd_r)$ are self-dual
$o(2m,\mbb{C})$-submodules for $r\in\ol{1,m-2}$. In particular, the
ternary weight code ${\cal C}_2$ of $o(2m,\mbb{C})$ on ${\cal A}_2$
is given by the weight matrix on its subspace
 $${\cal A}_{2,1}=\sum_{1\leq i<j\leq m}(\mbb{C}\sta_i\sta_j+\mbb{C}\sta_i\sta_{m+j}).\eqno(3.9)$$
We take any order
$$\{x_1,x_2,\cdots,x_{m(m-1)}\}=\{\sta_i\sta_j,\sta_i\sta_{m+j}\mid 1\leq i<j\leq m\}\eqno(3.10)$$
and write
$$(E_{i,i}-E_{m+i,m+i})(x_j)=c_{i,j}(2)x_j,\qquad C_2=(c_{i,j}(2))_{m\times m(m-1)}.\eqno(3.11)$$
Moreover,
$$\mbox{the weight matrix on}\; {\cal A}_2\;\mbox{is equivalent to}\;(C_2,-C_2). \eqno(3.12)$$
Since
$$\sum_{i=1}^m\mbb{F}_3h_i=\sum_{i=1}^m\mbb{F}_3(E_{i,i}-E_{m+i,m+i}),\eqno(3.13)$$
$C_2$ is a generator matrix of the ternary code ${\cal C}_2$. Denote
by $\zeta_i$ the $i$th row of $C_2$.  By (3.8) and (3.12), any
nonzero codeword in ${\cal C}_2$ has the same weight as the codeword
$$u(t)=\sum_{i=1}^t\zeta_t\qquad\mbox{for some}\;\;t\in\ol{1,m}.
\eqno(3.14)$$ Moreover,
$$f(t)=\wt u(t)={t\choose 2}+2t(m-t)=\frac{(4m-1)t-3t^2}{2}\eqno(3.15)$$ So $f(t)$ has only local maximum and it
attains the absolute  minimum at the boundary points. Note that
$$f(1)=2(m-1),\qquad f(m)=\frac{m(m-1)}{2}.\eqno(3.16)$$
Hence
$$\mbox{the minimal distance of}\;\;{\cal
C}_2\;\;\mbox{is}\;\;2(m-1)\;\;\mbox{if}\;m\geq 4.\eqno(3.17)$$ \pse

{\bf Theorem 3.1}. {\it When $m=3m_1+1$ for some positive integer
$m_1$, the ternary weight code ${\cal C}_2$ of $o(2m,\mbb{C})$ on
${\cal A}_2$ is an orthogonal $[m(m-1),m,2(m-1)]$-code.}

{\it Proof}. Note that for $i,j\in\ol{1,m}$ with $i\neq j$,
$$\zeta_i\cdot\zeta_i=f(1)=6m_1,\;\;(\zeta_i+\zeta_j)\cdot(\zeta_i+\zeta_j)=f(2)=1+4(m-2)=4m-7=12(m_1-1).\eqno(3.18)$$
Thus
$$\zeta_i\cdot \zeta_j=\frac{f(2)-2f(1)}{2}=-6.\eqno(3.19)$$
Hence  ${\cal C}_2$ is an orthogonal ternary code. $\qquad\Box$
\psp

In particular, ${\cal C}_2$ is an orthogonal ternary
$[12,4,6]$-code when $m_1=1$, $[42,7,12]$-code when $m_1=2$ and
$[90,10,18]$-code when $m_1=3$. It can be proved that ${\cal C}_2$
is also the weight code on the adjoint module of $o(2m,\mbb{C})$.

The ternary weight code ${\cal C}_3$ of  $o(2m,\mbb{C})$ on ${\cal
A}_3$ is given by the weight matrix on its subspace
 $${\cal A}_{3,1}=\sum_{1\leq i<j<l\leq m}\mbb{C}\sta_i\sta_j\sta_l+\sum_{1\leq i<j\leq m}\;\sum_{l=1}^m
 \mbb{C}\sta_i\sta_j\sta_{m+l}.\eqno(3.20)$$
We take any order
\begin{eqnarray*} & &\{y_1,y_2,\cdots,y_{{m\choose 3}+m{m\choose 2}}\}\\ &=&\{\sta_i\sta_j\sta_l,
\sta_r\sta_s\sta_{m+q}\mid 1\leq i<j<l\leq m;\;1\leq r<s\leq
m;\;q\in\ol{1,m}\}\hspace{2.4cm}(3.21)\end{eqnarray*} and write
$$(E_{i,i}-E_{m+i,m+i})(y_j)=c_{i,j}(3)y_j,\qquad C_3=(c_{i,j}(3))_{m\times\left({m\choose 3}+m{m\choose 2}\right)}.
\eqno(3.22)$$ Moreover,
$$\mbox{the weight matrix on}\; {\cal A}_3\;\mbox{is equivalent to}\;(C_3,-C_3).
\eqno(3.23)$$

Denote by $\eta_i$ the $i$th row of $C_3$.  By (3.8) and (3.23), any
nonzero codeword in ${\cal C}_3$ has the same weight as the codeword
$$u(t)=\sum_{i=1}^t\eta_t\qquad\mbox{for some}\;\;t\in\ol{1,m}.
\eqno(3.24)$$ Moreover,
\begin{eqnarray*}\hspace{2cm}g(t)&=&\wt u(t)=(2m-t){t\choose 2}+2t{m-t\choose 2}+t(m-t)^2
\\ &=&\frac{t(t-1)(2m-t)+2t(m-t)(2m-2t-1)}{2}\\ &=&\frac{t}{2}[3t^2+3(1-2m)t+4(m^2-m)].
\hspace{5.2cm}(3.25)\end{eqnarray*} Observe that
$$g'(t)=\frac{1}{2}[9t^2+6(1-2m)t+4(m^2-m)]=\frac{1}{2}[(3t+1-2m)^2-1].\eqno(3.26)$$
Thus
$$g'(t_0)=0\lra t_0=\frac{2(m-1)}{3},\;\frac{2m}{3}.\eqno(3.27)$$
Since $g'(0)=(m^2-m)/2\geq 0$, $t=2(m-1)/3$ is a point of local
maximum and $t=2m/3$ is a point of local minimum. We calculate
$$g(1)=(m-1)(2m-3),\qquad g(m)=\frac{m^2(m-1)}{2},\qquad
g(2m/3)=\frac{2}{9}m^2(2m-3).\eqno(3.28)$$ Note that $g(m)\geq g(1)$
and $g(2m/3)\geq g(1)$ if $m\geq 3$. \psp

{\bf Theorem 3.2}. {\it Let $m\geq 3$. The ternary weight code
${\cal C}_3$ of  $o(2m,\mbb{C})$ on ${\cal A}_3$ is of type
$[m(m-1)(2m-1)/3,m,(m-1)(2m-3)]$. Moreover, it is orthogonal if
 $m\not\equiv -1\;(\mbox{mod}\;3)$.}

{\it Proof}. Note
$$\eta_i\cdot\eta_i=g(1)=(m-1)(2m-3)\eqno(3.29)$$
and
$$(\eta_i+\eta_j)\cdot(\eta_i+\eta_j)=g(2)=2(2(m-2)^2+1)\eqno(3.30)$$
for $i,j\in\ol{1,m}$ such that $i\neq j$. Thus
$$\eta_i\cdot\eta_j=\frac{g(2)-2g(1)}{2}=3(2-m).\eqno(3.31)$$
So  ${\cal C}_3$ is orthogonal if $m\not\equiv -1\;(\mbox{mod}\;3).\qquad\Box$\psp

Remark that ${\cal C}_3$ is an orthogonal $[10,3,6]$-code when
$m=3$, $[28,4,15]$-code when $m=4$, $[110,6,45]$-code when $m=6$
and $[182,7,66]$-code when $m=7$.

Let ${\cal B}$ be the subalgebra of ${\cal A}$ generated by
$\{1_{\cal A},\sta_i\mid i\in\ol{1,m}\}$ and
$${\cal B}_r={\cal A}_r\bigcap {\cal B}\qquad\for\;\;r\in\ol{0,m}.\eqno(3.32)$$
The spin representation of $so(2m,\mbb{C})$ is given by
$$E_{i,j}-E_{m+j,m+i}=\sta_i\ptl_{\sta_j}-\frac{\dlt_{i,j}}{2}\qquad
\for\;\;i,j\in\ol{1,m},\eqno(3.33)$$
$$E_{m+s,r}-E_{m+r,s}=\ptl_{\sta_s}\ptl_{\sta_r},\qquad
 E_{r,m+s}-E_{s,m+r}=\sta_r\sta_s\eqno(3.34)$$
 for $1\leq r<s\leq m$. Then the subspace
$${\cal V}=\sum_{i=1}^{\llbracket m/2\rrbracket}{\cal B}_{m-i}\eqno(3.35)$$
is the irreducible module with highest weight $\lmd_m$, that is,
${\cal V}\cong V_{D_m}(\lmd_m)$.

If $m=2m_1+1$ is odd, then
$$\{\sta_{i_1}\cdots\sta_{i_{m-2r}}\mid r\in
\ol{0,m_1};\;1\leq i_1<\cdots< i_{m-2r}\leq m\}\eqno(3.36)$$ forms a
weight-vector basis of ${\cal V}$. When $m=2m_1$ is even,
$$\{1,\sta_{i_1}\cdots\sta_{i_{m-2r}}\mid r\in
\ol{0,m_1-1};\;1\leq i_1<\cdots< i_{m-2r}\leq m\}\eqno(3.37)$$ is
a weight-vector basis of ${\cal V}$. Take any order
$\{z_1,z_2,...,z_{2^{m-1}}\}$ of the above base vectors. Denote
$$(E_{r,r}-E_{m+r,m+r})(z_i)=q_{r,i}z_i,\qquad
C({\cal V})=(q_{r,i})_{m\times 2^{m-1}}.\eqno(3.38)$$ Note that
$$\frac{1}{2}\equiv -1\qquad\mbox{in}\;\;\mbb{F}_3.\eqno(3.39)$$

Denote by $\xi_r$ the $r$th row of the weight matrix $C({\cal V})$.
Set
$$\bar u=\sum_{r=1}^{m-1}\xi_r-\xi_m,\;\;
u(t)=\sum_{i=1}^t\xi_i\qquad\for\;\;t\in\ol{1,m}.\eqno(3.40)$$ Then
any  nonzero codeword in ${\cal C}_3({\cal V})$ is conjugated to
some $u(t)$ or $\bar u$ under the action of the Weyl group of
$o(2m,\mbb{C})$ (cf. (1.10) and (1.11)). It has the same weight as
$u(t)$ or $\bar u$. We calculate
$$\wt u(1)=2^{m-1},\qquad\wt u(2)=2^{m-2}.\eqno(3.41)$$
Moreover, we have the following more general estimates. For any
positive integer $k>2$, we always have
$${k\choose l-1}+{k\choose l+1}>{k\choose
l}\qquad\for\;\;l\in\ol{0,k},\eqno(3.42)$$ where we treat ${k\choose
-1}={k\choose k+1}=0$. If $t=3t_1$ for some positive integer $t_1$,
we have
\begin{eqnarray*} \wt
u(t)&=&2^{m-3t_1-1}\sum_{i=0}^{t_1}\left[{3t_1\choose
6i+1}+{3t_1\choose 6i+2}+{3t_1\choose 6i+4}+{3t_1\choose 6i+5}\right]\\
&>& 2^{m-3t_1-1}\sum_{i=0}^{t_1}\left[{3t_1\choose
6i+1}+{3t_1\choose 6i+3}+{3t_1\choose 6i+5}\right]=
2^{m-2}.\hspace{2.3cm}(3.43)\end{eqnarray*} When $t=3t_1+1$ for some
positive integer $t_1$, we obtain
\begin{eqnarray*} \wt
u(t)&=&2^{m-3t_1-2}\sum_{i=0}^{t_1}\left[{3t_1+1\choose
6i}+{3t_1+1\choose 6i+1}+{3t_1+1\choose 6i+3}+{3t_1+1\choose 6i+4}\right]\\
&>& 2^{m-3t_1-2}\sum_{i=0}^{t_1}\left[{3t_1+1\choose
6i}+{3t_1+1\choose 6i+2}+{3t_1+1\choose 6i+4}\right]=
2^{m-2}.\hspace{1.7cm}(3.44)\end{eqnarray*} If $t=3t_1+2$ for some
positive integer $t_1$, we get
\begin{eqnarray*} \wt
u(t)&=&2^{m-3t_1-3}\sum_{i=0}^{t_1}\left[{3t_1+2\choose
6i}+{3t_1+2\choose 6i+2}+{3t_1+2\choose 6i+3}+{3t_1+2\choose 6i+5}\right]\\
&>& 2^{m-3t_1-3}\sum_{i=0}^{t_1}\left[{3t_1+2\choose
6i}+{3t_1+2\choose 6i+2}+{3t_1+2\choose 6i+4}\right]=
2^{m-2}.\hspace{1.7cm}(3.45)\end{eqnarray*}

Let $k$ be  a positive integer. We have
$${2k\choose i}+{2k\choose i+4}>{2k\choose i+1}\eqno(3.46)$$
if $i\leq k-3$ or $i\geq k$. Moreover,
$${2k\choose k-2}+{2k\choose k+2}-{2k\choose k-1}=\frac{k-4}{k-1}{2k\choose
k-2},\eqno(3.47)$$
$${2k\choose k-1}+{2k\choose k+3}-{2k\choose
k}=\frac{k^3-4k^2-3k-6}{k(k-1)(k-2)}{2k\choose k-3}.\eqno(3.48)$$
Thus (3.46) always holds if $k\geq 5$. Furthermore,
$${2k+1\choose i}+{2k+1\choose i+4}>{2k+1\choose i+1}\eqno(3.49)$$
if $i\neq  k-1$. Observe that
$${2k+1\choose k-1}+{2k+1\choose i+3}-{2k+1\choose
k}=\frac{k^2-3k-6}{k(k-1)}{2k+1\choose k-2}.\eqno(3.50)$$ So (3.49)
holds whenever $k\geq 5$. Therefore,
$${k\choose i}+{k\choose i+4}>{k\choose i+1}\qquad\mbox{if}\;\;k\geq 10.\eqno(3.51)$$

If $m=3m_1$ for some positive integer $m_1$,
\begin{eqnarray*} \wt\bar u&=&\sum_{i=0}^m\left[{m\choose 6i}+{m-1\choose
6i+1}+{m-1\choose 6i+4}\right]\\ &=&\sum_{i=0}^m\left[{m-1\choose
6i}+{m-1\choose 6i+5}+{m-1\choose 6i+1}+{m-1\choose
6i+4}\right],\hspace{3.4cm}(3.52)\end{eqnarray*} which is $>
2^{m-2}$ if $m_1\geq 4$ by (3.51). When $m=3m_1+1$ for some positive
integer $m_1$,
\begin{eqnarray*} \wt\bar u&=&\sum_{i=0}^m\left[{m-1\choose
6i}+{m\choose 6i+2}+{m-1\choose 6i+3}\right]\\
&=&\sum_{i=0}^m\left[{m-1\choose 6i}+{m-1\choose 6i+1}+{m-1\choose
6i+2}+{m-1\choose 6i+3}\right]\\
&=&1+\sum_{i=0}^m\left[{m-1\choose 6i+1}+{m-1\choose
6i+3}+{m-1\choose 6i+2}+{m-1\choose
6i+6}\right],\hspace{2.7cm}(3.53)\end{eqnarray*} which is again
$>2^{m-2}$ if $m_1\geq 4$ by (3.51). Assuming  $m=3m_1+2$ for some
positive integer $m_1$, we have
\begin{eqnarray*} \wt\bar u&=&\sum_{i=0}^m\left[{m-1\choose
6i+2}+{m\choose 6i+4}+{m-1\choose 6i+5}\right]\\ &=&{m-1\choose
3}+\sum_{i=0}^m\left[{m-1\choose 6i+2}+{m-1\choose 6i+4}+{m-1\choose
6i+5}+{m-1\choose 6i+9}\right] ,\hspace{1.2cm}(3.54)\end{eqnarray*}
which is $> 2^{m-2}$ if $m_1\geq 3$ by (3.51).
 Moreover, we have the
following table:
\begin{center}{\bf Table 3.1}\end{center}
\begin{center}\begin{tabular}{|c|c|c|c|c|c|c|c|}\hline
m&4&5&6&7&8&9&10
\\\hline $\wt\bar u$&8&11&12&43&112&171&260\\\hline\end{tabular}\end{center}
In summary, we have:\psp

{\bf Theorem 3.3}. {\it Let $m>3$ be an integer. The ternary code
${\cal C}_3({\cal V})$ is of type $[2^{m-1},m,2^{m-2}]$ if $m\neq 6$
and of type $[32,6,12]$ when $n=6$}.\psp

We remark that the spin module ${\cal V}$ is self-dual if and only if $m$ is even. \psp

{\bf Corollary 3.4}. {\it When $m=6m_1+2$ for some positive
integer $m_1$, the ternary weight code of $o(2m,\mbb{C})$ on
$o(2m,\mbb{C})+{\cal V}$ is an orthogonal ternary
$[m(m-1)+2^{m-2},m,4m-7+2^{m-3}]$-code. If
 $m=6m_1+3$ for some positive integer $m_1$, the ternary weight code of
 $o(2m,\mbb{C})$ on
$o(2m,\mbb{C})+{\cal V}$ is an orthogonal ternary
$[2m(m-1)+2^{m-1},m,8m-14+2^{m-2}]$-code. In the case $m=6m_1+5$
and $m=6m_1+12$ for some nonnegative integer $m_1$, the code
${\cal C}_2\oplus {\cal C}_3({\cal V})$ is an orthogonal ternary
$[m(m-1)+2^{m-1},m, 4m-7+2^{m-2}]$-code. When $m=6$, the code
${\cal C}_2\oplus {\cal C}_3({\cal V})$ is an orthogonal ternary
$[62,6, 27]$-code.}

{\it Proof}. Suppose $m=6m_1+2$ for some positive integer $m_1$.
Then the weight matrix of $o(2m,\mbb{C})$ on $o(2m,\mbb{C})+{\cal
V}$ is equivalent to $(A,-A)$, where $A$ generates the weight code
${\cal C}$ of $o(2m,\mbb{C})+{\cal V}$. Moreover, ${\cal C}$ is
orthogonal if and only if the matrix $(A,-A)$ generates an
orthogonal code. But
$$(A,-A)\sim (C_2,C_2,C({\cal V})).\eqno(3.55)$$
Note that
$$\wt (\zeta_i,\zeta_i,\xi_i)=2f(1)+2^{m-1}=4(m-1)+2^{m-1}\equiv 1+(-1)^{6m_1+1}\equiv 0
\;\;(\mbox{mod}\;3),\eqno(3.56)$$
\begin{eqnarray*} & &\wt (\zeta_i+\zeta_j,\zeta_i+\zeta_j,\xi_i+\xi_j)
\\ &=&2f(2)+2^{m-2}=8m-14+2^{m-2}\equiv 2+(-1)^{6m_1}\equiv 0
\;\;(\mbox{mod}\;3)\hspace{2.9cm}(3.57)\end{eqnarray*} for
$i,j\in\ol{1,m}$ with $i\neq j$ by (3.16) and (3.41). Thus
$$(\zeta_i,\zeta_i,\xi_i)\cdot(\zeta_i,\zeta_i,\xi_i)\equiv \wt (\zeta_i,\zeta_i,\xi_i)\equiv 0,\eqno(3.58)$$
\begin{eqnarray*} & &(\zeta_i,\zeta_i,\xi_i)\cdot(\zeta_j,\zeta_j,\xi_j)
\\ &\equiv& -[\wt (\zeta_i+\zeta_j,\zeta_i+\zeta_j,\xi_i+\xi_j)-\wt (\zeta_i,\zeta_i,\xi_i)-
\wt (\zeta_j,\zeta_j,\xi_j)]\equiv
0\hspace{2.5cm}(3.59)\end{eqnarray*} by (3.39). Thus ${\cal C}$ is
orthogonal. Note
$$f(2)=4m-7\leq \frac{m(m-1)}{2}=f(m)\qquad\mbox{if}\;\;m\geq 7.\eqno(3.60)$$
Thus
$$f(2)\leq f(t)\qquad\for\;\;t\in\ol{2,m}.\eqno(3.61)$$
By (3.8),
$$\wt (\sum_{i=1}^{m-1}\zeta_i-\zeta_m)=f(m)\geq f(2).\eqno(3.62)$$
Thus the minimum distance of ${\cal C}$ is
$$\min\{f(1)+2^{m-2},f(2)+2^{m-3}\}=4m-7+2^{m-3}\qquad\mbox{if}\;\;m\geq 6.\eqno(3.63)$$
This proves the first conclusion. The other conclusions for $m\geq 7$ can be proved similarly.

In the case $m=5$, we have
\begin{center}{\bf Table 3.2}\end{center}
\begin{center}\begin{tabular}{|c|c|c|c|c|c|}\hline
t&1&2&3&4&5
\\\hline f(t)&8&13&15&14&10\\\hline\end{tabular}\end{center}
and on the ${\cal V}$,
\begin{center}{\bf Table 3.3}\end{center}
\begin{center}\begin{tabular}{|c|c|c|c|c|c|}\hline
t&1&2&3&4&5
\\\hline \wt u(t)&16&8&12&10&11\\\hline\end{tabular}\end{center}
By Tables 3.1-3.3 and the fact $\wt
(\sum_{i=1}^4\zeta_i-\zeta_5)=f(5)$ in ${\cal C}_3({\cal A}_2)$,
the third conclusion holds for $m=5$.

If $m=6$,
\begin{center}{\bf Table 3.4}\end{center}
\begin{center}\begin{tabular}{|c|c|c|c|c|c|c|}\hline
t&1&2&3&4&5&6
\\\hline f(t)&10&17&21&22&20&15\\\hline\end{tabular}\end{center}
and on the ${\cal V}$,
\begin{center}{\bf Table 3.5}\end{center}
\begin{center}\begin{tabular}{|c|c|c|c|c|c|c|}\hline
t&1&2&3&4&5&6
\\\hline \wt u(t)&32&16&24&20&22&21\\\hline\end{tabular}\end{center}
By Tables 3.1, 3.4, and 3.5,  and the fact $\wt
(\sum_{i=1}^5\zeta_i-\zeta_6)=f(6)$ in ${\cal C}_3({\cal A}_2)$,
the last conclusion holds. $\qquad\Box$ \psp

When $m=8$, the ternary weight code of  $o(16,\mbb{C})$ on
$o(16,\mbb{C})+{\cal V}$ is a ternary orthogonal  $[120,8,57]$-code,
which will later be proved also to be the ternary weight code of
$E_8$ on its adjoint module. If $m=9$, the ternary weight code of
$o(18,\mbb{C})$ on $o(18,\mbb{C})+{\cal V}$ is a ternary orthogonal
$[400,8,186]$-code. When $m=5$, the code ${\cal C}_2\oplus {\cal
C}_3({\cal V})$ is a  ternary orthogonal $[36,5,21]$-code, which
will later be proved also to be the ternary weight code of $E_6$ on
its adjoint module. In the case $m=11$, the code ${\cal C}_2\oplus
{\cal C}_3({\cal V})$ is a ternary orthogonal $[1134,8,549]$-code.

\section{Representations of $F_4$ and Ternary Codes}

In this section, we study the ternary weight codes of $F_4$ on its
minimal irreducible module and adjoint module.

We go back to the settings in (2.2)-(2.4) with $n=4$.  The root
system of $F_4$ is
$$\Phi_{F_4}=\left\{\pm \ves_i,\pm \ves_i\pm \ves_j,\frac{1}{2}(\pm \ves_1\pm \ves_2\pm
\ves_3\pm \ves_4)\mid i\neq j\right\}\eqno(4.1)$$  and the
positive simple roots are
$$\al_1=\ves_2-\ves_3,\al_2=\ves_3-\ves_4,\al_3=\ves_4,\al_4=\frac{1}{2}(\ves_1-\ves_2-\ves_3-\ves_4).
\eqno(4.2)$$ The corresponding  Dynkin diagram  is

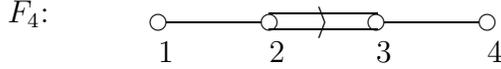
\begin{picture}(60,12)\put(2,0){$F_4$:}
\put(21,0){\circle{2}}\put(21,-5){1}\put(22,0){\line(1,0){12}}
\put(35,0){\circle{2}}\put(35,-5){2}\put(35,1.2){\line(1,0){13.6}}
\put(35,-0.8){\line(1,0){13.6}}\put(41,-1){$\ra$}\put(48.5,0){\circle{2}}\put(48.5,-5){3}\put(49.5,0)
{\line(1,0){12}}\put(62.5,0){\circle{2}}\put(62.5,-5){4}
\end{picture}
\vspace{0.6cm}

The Weyl group ${\cal W}_{F_4}$ of $F_4$ contains the permutation
group $S_4$ on the sub-indices of $\ves_i$ and all reflections with
respect to the coordinate hyperplanes. Moreover, there is an
identification:
$$h_1\leftrightarrow\al_1,\;h_2\leftrightarrow\al_2,\;h_3\leftrightarrow 2\al_3,\;h_4\leftrightarrow 2\al_4\eqno(4.3)$$
(e.g, cf. [7]). Thus
$${\cal H}_2=\sum_{i=1}^4\mbb{F}_2h_i=\sum_{i=1}^4\mbb{F}_2\ves_i.\eqno(4.4)$$
Moreover,
$${\cal H}_2=\{{\cal W}_{F_4}(h_1),{\cal W}_{F_4}(h_1+h_3),{\cal W}_{F_4}(h_3),{\cal W}_{F_4}(h_4)\}.\eqno(4.5)$$
The basic (minimal) irreducible module $V_{F_4}$ of the 52-dimensional Lie
algebra ${\cal G}^{F_4}$ has a basis $\{x_i\mid\ol{1,26}\}$ and with
the representation determined by the following formulas in terms of
differential operators:
$$E_{\al_1}|_V=x_4\ptl_{x_6}+x_5\ptl_{x_8}+x_7\ptl_{x_9}-
x_{18}\ptl_{x_{20}}-x_{19}\ptl_{x_{22}}-x_{21}\ptl_{x_{23}},\eqno(4.6)$$
$$E_{\al_2}|_V=x_3\ptl_{x_4}+x_8\ptl_{x_{10}}+x_9\ptl_{x_{11}}
-x_{16}\ptl_{x_{18}}-x_{17}\ptl_{x_{19}}-x_{23}\ptl_{x_{24}},\eqno(4.7)$$
\begin{eqnarray*}\hspace{1cm}E_{\al_3}|_V&=&-x_2\ptl_{x_3}-x_4\ptl_{x_5}-x_6\ptl_{x_8}
+x_{10}\ptl_{x_{12}}+x_{11}(\ptl_{x_{13}}-2\ptl_{x_{14}})\\
&
&-x_{14}\ptl_{x_{16}}-x_{15}\ptl_{x_{17}}+x_{19}\ptl_{x_{21}}+x_{22}\ptl_{x_{23}}+x_{24}\ptl_{x_{25}},
\hspace{3.3cm}(4.8)\end{eqnarray*}
\begin{eqnarray*}\hspace{1cm}E_{\al_4}|_V&=&-x_1\ptl_{x_2}-x_5\ptl_{x_7}-x_8\ptl_{x_9}-
x_{10}\ptl_{x_{11}}+x_{12}(\ptl_{x_{14}}-2\ptl_{x_{13}})
\\ &
&-x_{13}\ptl_{x_{15}}+x_{16}\ptl_{x_{17}}+x_{18}\ptl_{x_{19}}+x_{20}\ptl_{x_{22}}
+x_{25}\ptl_{x_{26}},\hspace{3.3cm}(4.9)\end{eqnarray*}
$$E_{-\al_1}|_V=-x_6\ptl_{x_4}-x_8\ptl_{x_5}-x_9\ptl_{x_7}+
x_{20}\ptl_{x_{18}}+x_{22}\ptl_{x_{19}}+x_{23}\ptl_{x_{21}},\eqno(4.10)$$
$$E_{-\al_2}|_V=-x_4\ptl_{x_3}-x_{10}\ptl_{x_8}-x_{11}\ptl_{x_{9}}
+x_{18}\ptl_{x_{16}}+x_{19}\ptl_{x_{17}}+x_{24}\ptl_{x_{23}},\eqno(4.11)$$
\begin{eqnarray*}\hspace{1cm}E_{-\al_3}|_V&=&x_3\ptl_{x_2}+x_5\ptl_{x_4}+x_8\ptl_{x_6}
-x_{12}\ptl_{x_{10}}+x_{16}(2\ptl_{x_{14}}-\ptl_{x_{13}})\\
&&+x_{14}\ptl_{x_{11}}+x_{17}\ptl_{x_{15}}-x_{21}\ptl_{x_{19}}-x_{23}\ptl_{x_{22}}-x_{25}\ptl_{x_{24}},
\hspace{2.9cm}(4.12)\end{eqnarray*}
\begin{eqnarray*}\hspace{1cm}E_{-\al_4}|_V&=&x_2\ptl_{x_1}+x_7\ptl_{x_5}+x_9\ptl_{x_8}+
x_{11}\ptl_{x_{10}}+x_{15}(2\ptl_{x_{13}}-\ptl_{x_{14}})
\\ &
&+x_{13}\ptl_{x_{12}}-x_{17}\ptl_{x_{16}}-x_{19}\ptl_{x_{18}}-x_{22}\ptl_{x_{20}}
-x_{26}\ptl_{x_{25}},\hspace{2.8cm}(4.13)\end{eqnarray*}
\begin{eqnarray*}\hspace{1cm}h_1|_V&=&x_4\ptl_{x_4}+x_5\ptl_{x_5}-x_6\ptl_{x_6}+x_7\ptl_{x_7}-x_8\ptl_{x_8}-x_9\ptl_{x_9}
+x_{18}\ptl_{x_{18}}\\
&&+x_{19}\ptl_{x_{19}}-x_{20}\ptl_{x_{20}}+x_{21}\ptl_{x_{21}}-x_{22}\ptl_{x_{22}}
-x_{23}\ptl_{x_{23}},\hspace{3.4cm}(4.14)
\end{eqnarray*}
\begin{eqnarray*}\hspace{1cm}h_2|_V&=&x_3\ptl_{x_3}-x_4\ptl_{x_4}+x_8\ptl_{x_8}
+x_9\ptl_{x_9}-x_{10}\ptl_{x_{10}}-x_{11}\ptl_{x_{11}}
+x_{16}\ptl_{x_{16}}\\&
&+x_{17}\ptl_{x_{17}}-x_{18}\ptl_{x_{18}}-x_{19}\ptl_{x_{19}}+x_{23}\ptl_{x_{23}}
-x_{24}\ptl_{x_{24}},\hspace{3.4cm}(4.15)
\end{eqnarray*}
\begin{eqnarray*}\hspace{1cm}h_3|_V&=&x_2\ptl_{x_2}-x_3\ptl_{x_3}+x_4\ptl_{x_4}-x_5\ptl_{x_5}+x_6\ptl_{x_6}-x_8\ptl_{x_8}
+x_{10}\ptl_{x_{10}}\\&
&+2x_{11}\ptl_{x_{11}}-x_{12}\ptl_{x_{12}}+x_{15}\ptl_{x_{15}}-2x_{16}\ptl_{x_{16}}-x_{17}\ptl_{x_{17}}
+x_{19}\ptl_{x_{19}}
\\
&&-x_{21}\ptl_{x_{21}}+x_{22}\ptl_{x_{22}}-x_{23}\ptl_{x_{23}}+x_{24}\ptl_{x_{24}}
-x_{25}\ptl_{x_{25}},\hspace{3.4cm}(4.16)
\end{eqnarray*}
\begin{eqnarray*}\hspace{1cm}h_4|_V&=&x_1\ptl_{x_1}-x_2\ptl_{x_2}+x_5\ptl_{x_5}-x_7\ptl_{x_7}+x_8\ptl_{x_8}-x_9\ptl_{x_9}
+x_{10}\ptl_{x_{10}}\\&
&-x_{11}\ptl_{x_{11}}+2x_{12}\ptl_{x_{12}}-2x_{15}\ptl_{x_{15}}+x_{16}\ptl_{x_{16}}-x_{17}\ptl_{x_{17}}
+x_{18}\ptl_{x_{18}}
\\
&&-x_{19}\ptl_{x_{19}}+x_{20}\ptl_{x_{20}}-x_{22}\ptl_{x_{22}}+x_{25}\ptl_{x_{25}}
-x_{26}\ptl_{x_{26}}\hspace{3.4cm}(4.17)
\end{eqnarray*}
(e.g., cf. [26])

The module $V_{F_4}$ is self-dual. The weight matrix of $V_{F_4}$ is
$(A_{F_4},-A_{F_4})$ with
$$A_{F_4}=\left[\begin{array}{rrrrrrrrrrrr}0&0&0&1&1&-1&1&-1&-1&0&0&0\\
0&0&1&-1&0&0&0&1&1&-1&-1&0\\
0&1&-1&1&-1&1&0&-1&0&1&2&-1\\
1&-1&0&0&1&0&-1&1&-1&1&-1&2\end{array}\right].\eqno(4.18)$$ \pse

{\bf Theorem 4.1}. {\it The ternary weight code ${\cal C}_{F_4,1}$
(generated by $A_{F_4}$) of $F_4$ on $V_{F_4}$ is an orthogonal
[12,4,6]-code}.

 {\it Proof}. Denote by $\xi_i$ the $i$th row of the matrix $A_{F_4}$. Then
$$\wt\xi_1=6,\qquad \wt (\xi_1+\xi_3)=\wt\xi_3=\wt\xi_4=9.\eqno(4.19)$$
According to (4.5), any nonzero codeword in ${\cal C}_{F_4,1}$ has
weight 6 or 9. By an argument as (3.29)-(3.31), ${\cal C}_{F_4,1}$
is orthogonal.$\qquad\Box$\psp

Next we consider the adjoint representation of $F_4$. Its weight
code ${\cal C}_{F_4,2}$ is determined by the set $\Phi_{F_4}^+$ of positive roots. The followings are positive
roots of $F_4$:
$$\al_1,\;\;\al_2,\;\;\al_3,\;\;\al_4,\;\;\al_1+\al_2,\;\;\al_2+\al_3,\;\;\al_3+\al_4,\;\;
\al_1+\al_2+\al_3,\;\;\al_2+\al_3+\al_4,\eqno(4.20)$$
$$\al_2+2\al_3,\;\;\al_1+\al_2+2\al_3,\;\;\al_2+2\al_3+\al_4,\;\;\al_1+\al_2+\al_3+\al_4,
\eqno(4.21)$$
$$\al_1+2\al_2+2\al_3,\;\;\al_1+\al_2+2\al_3+\al_4,\;\;\al_2+2\al_3+2\al_4,\;\;\al_1+2\al_2+2\al_3+\al_4,\eqno(4.22)$$
$$\al_1+\al_2+2\al_3+2\al_4,\;\;\al_1+2\al_2+2\al_3+2\al_4,\;\;\al_1+2\al_2+3\al_3+\al_4,\;\;\al_1+2\al_2+3\al_3+2\al_4,
\eqno(2.22)$$
$$\al_1+2\al_2+4\al_3+2\al_4,\;\;\al_1+3\al_2+4\al_3+2\al_4,\;\;2\al_1+3\al_2+4\al_3+2\al_4.
\eqno(4.23)$$ Let $E_\al$ be a root vector associated with the root
$\al$. The weight matrix $B_{F_4}$ on
$\sum_{\al\in\Phi_{F_4}^+}\mbb{F}E_\al$ is given by
$${\tiny \left[\begin{array}{rrrrrrrrrrrrrrrrrrrrrrrr}2&$-1$&0&0&1&$-1$&0&1&$-1$&$-1$&1&$-1$&1&0&1&$-1$&0&1&0&0&0&0&$-1$&1\\
$-1$&2&$-1$&0&1&1&$-1$&0&1&0&$-1$&0&0&1&$-1$&0&1&$-1$&1&0&0&$-1$&1&0\\
0&$-2$&2&$-1$&$-2$&0&1&0&$-1$&2&2&1&$-1$&0&1&0&$-1$&0&$-2$&1&0&2&0&0\\
0&0&$-1$&2&0&$-1$&1&$-1$&1&$-2$&$-2$&0&1&$-2$&0&2&0&2&2&$-1$&1&0&0&0\end{array}\right].}\eqno(4.24)$$
\pse

{\bf Theorem 4.2}. {\it The ternary weight code ${\cal C}_{F_4,2}$
(generated by $B_{F_4}$) of $F_4$ on its adjoint module  is an
orthogonal $[24,4,15]$-code.}\psp

{\it Proof}. Denote by $\eta_i$ the $i$th row of the above matrix. Then
$$\wt \eta_i=15,\qquad\wt(\eta_1+\eta_3)=18.\eqno(4.25)$$
According to (4.5), any nonzero codeword in ${\cal C}_{F_4,2}$ has
weight 15 or 18. By an argument as (3.29)-(3.31) , ${\cal
C}_{F_4,2}$ is orthogonal.$\qquad\Box$\psp

\section{Representations of $E_6$ and Ternary Codes}

In this section, we investigate  the ternary weight codes of $E_6$
on its minimal irreducible module and adjoint module.

First we give a lattice construction of the exceptional simple Lie
algebras of type $E$. Let $\{\al_i\mid i\in\ol{1,m}\}$ be the simple
positive roots of type $E_m$.
 Set
$$Q_{E_m}=\sum_{i=1}^m\mbb{Z}\al_i,\eqno(5.1)$$ the root lattice of type
$E_m$. Denote by $(\cdot,\cdot)$ the symmetric $\mbb{Z}$-bilinear
form on $Q_{E_m}$ such that the root system
$$\Phi_{E_m}=\{\al\in Q_{E_m}\mid (\al,\al)=2\}.\eqno(5.2)$$
Define $F(\cdot,\cdot):\; Q_{E_m}\times  Q_{E_m}\rta \{\pm 1\}$ by
$$F(\sum_{i=1}^mk_i\al_i,\sum_{j=1}^ml_j\al_j)=(-1)^{\sum_{i=1}^mk_il_i+\sum_{m\geq i>j\geq 1}k_il_j
(\al_i,\al_j)},\qquad k_i,l_j\in\mbb{Z}.\eqno(5.3)$$
Denote
$$H_{E_m}=\sum_{i=1}^m\mbb{C}\al_i.\eqno(5.4)$$
The simple Lie algebra of type $E_m$ is
$${\cal
G}^{E_m}=H_{E_m}\oplus\bigoplus_{\al\in
\Phi_{E_m}}\mbb{C}E_{\al}\eqno(5.5)$$ with the Lie bracket
$[\cdot,\cdot]$ determined by:
 $$[H_{E_m},H_{E_m}]=0,\;\;[h,E_{\al}]=(h,\al)E_{\al},\;\;[E_{\al},E_{-\al}]=-\al,
 \eqno(5.6)$$
 $$[E_{\al},E_{\be}]=\left\{\begin{array}{ll}0&\mbox{if}\;\al+\be\not\in \Phi_{E_m},\\
 F(\al,\be)E_{\al+\be}&\mbox{if}\;\al+\be\in\Phi_{E_m}.\end{array}\right.\eqno(5.7)$$
for $\al,\be\in\Phi_{E_m}$ and $h\in H_{E_m}$ (e.g., cf. [8], [25]).
Moreover,
$$h_i=\al_i\qquad\for\;\;i\in\ol{1,m}.\eqno(5.8)$$

Recall the settings in (2.2)-(2.4). Taking $n=7$, we have the
following root system of $E_6$:
$$\Phi_{E_6}=\left\{\ves_i-\ves_j,\frac{1}{2}(\sum_{s=1}^6\iota_s\ves_s\pm
\sqrt{2}\ves_7),\pm \sqrt{2}\ves_7\mid i,j\in\ol{1,6},i\neq
j;\iota_s=\pm 1;\sum_{i=1}^6\iota_i=0\right\}\eqno(5.9)$$
and the simple positive roots are
$$\al_1=\ves_1-\ves_2,\;\;\al_2=\frac{1}{2}(\sum_{j=1}^3(\ves_{3+j}-\ves_j)+\sqrt{2}\ves_7),\;\;
\al_i=\ves_{i-1}-\ves_i,\qquad i\in\ol{3,6}.\eqno(5.10)$$
The Dynkin diagram
is:\\
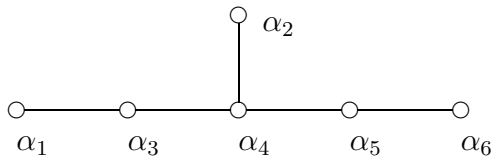
\begin{picture}(80,18)
\put(21,0){\circle{2}}\put(21,
-5){$\al_1$}\put(22,0){\line(1,0){12}}\put(35,0){\circle{2}}\put(35,
-5){$\al_3$}\put(36,0){\line(1,0){12}}\put(49,0){\circle{2}}\put(49,
-5){$\al_4$}\put(49,1){\line(0,1){10}}\put(49,12){\circle{2}}\put(52,10){$\al_2$}\put(50,0){\line(1,0){12}}
\put(63,0){\circle{2}}\put(63,-5){$\al_5$}\put(64,0){\line(1,0){12}}\put(77,0){\circle{2}}\put(77,
-5){$\al_6$}
\end{picture}\vspace{0.8cm}

Note
$${\cal H}_{E_6,3}=\sum_{i=1}^6\mbb{F}_3h_i=\{\sum_{i=1}^6\iota_i\ves_i+\iota_7\sqrt{2}\ves_7\mid
\iota_r\in\mbb{F}_3,\;\sum_{i=1}^6\iota_i=0\}.\eqno(5.11)$$
Moreover, the Weyl group ${\cal W}_{E_6}$ contains the permutation
group $S_6$ on the first six sub-indices  of $\ves_i$ and the
reflection
$$\sum_{i=1}^6\iota_i\ves_i+\iota_7\sqrt{2}\ves_7\mapsto \sum_{i=1}^6\iota_i\ves_i-\iota_7\sqrt{2}\ves_7.\eqno(5.12)$$
So
$${\cal H}_{E_6,3}={\cal W}_{E_6}(\{\sum_{i=1}^s\ves_i-\sum_{j=1}^t\ves_{s+j}+\iota\sqrt{2}\ves_7,\sqrt{2}\ves_7\mid
\iota=0,1;s-t\equiv 0\;(\mbox{mod}\;3)\}).\eqno(5.13)$$ The
27-dimensional basic irreducible module $V_{E_6}$ of weight $\lmd_1$
for $E_6$ has a basis $\{x_i\mid i\in\ol{1,27}\}$ with the
representation formulas determined by
$$E_{\al_1}|_V=-x_1\ptl_{x_2}+x_{11}\ptl_{x_{14}}+x_{15}\ptl_{x_{17}}
+x_{16}\ptl_{x_{19}}+x_{18}\ptl_{x_{21}}+x_{20}\ptl_{x_{23}},\eqno(5.14)$$
$$E_{\al_2}|_V=-x_4\ptl_{x_6}-x_5\ptl_{x_7}-x_8\ptl_{x_{10}}
+x_{18}\ptl_{x_{20}}+x_{21}\ptl_{x_{23}}+x_{22}\ptl_{x_{24}},\eqno(5.15)$$
$$E_{\al_3}|_V=-x_2\ptl_{x_3}+x_9\ptl_{x_{11}}
+x_{12}\ptl_{x_{15}}+x_{13}\ptl_{x_{16}}+x_{21}\ptl_{x_{22}}+x_{23}
\ptl_{x_{24}},\eqno(5.16)$$
$$E_{\al_4}|_V=-x_3\ptl_{x_4}-x_7\ptl_{x_9}-x_{10}\ptl_{x_{12}}
-x_{16}\ptl_{x_{18}}-x_{19}\ptl_{x_{21}}+x_{24}\ptl_{x_{25}},\eqno(5.17)$$
$$E_{\al_5}|_V=-x_4\ptl_{x_5}-x_6\ptl_{x_7}-x_{12}\ptl_{x_{13}}
-x_{15}\ptl_{x_{16}}
-x_{17}\ptl_{x_{19}}+x_{25}\ptl_{x_{26}},\eqno(5.18)$$
$$E_{\al_6}|_V=-x_5\ptl_{x_8}-x_7\ptl_{x_{10}}-x_9\ptl_{x_{12}}-x_{11}
\ptl_{x_{15}}
-x_{14}\ptl_{x_{17}}+x_{26}\ptl_{x_{27}},\eqno(5.19)$$
$$h_r|_{V_{E_6}}=\sum_{i=1}^{27}a_{r,i}x_i\ptl_{x_i}\eqno(5.20)$$
with $a_{r,i}$ given by the following table

\begin{center}{\bf \large Table 5.1}\end{center}
\begin{center}\begin{tabular}{|r||r|r|r|r|r|r||r||r|r|r|r|r|r|}\hline
$i$&$a_{1,i}$&$a_{2,i}$&$a_{3,i}$&$a_{4,i}$&$a_{5,i}$&$a_{6,i}$&$i$&$a_{1,i}$&$a_{2,i}$&$a_{3,i}$&$a_{4,i}$&$a_{5,i}$&
$a_{6,i}$
\\\hline\hline 1&1&0&0&0&0&0&2&$-1$&0&1&0&0&0\\\hline 3&0&0&$-1$&1&0&0&4&$0$&1&0&$-1$&1&0
\\\hline 5&0&1&$0$&0&$-1$&1&6&$0$&$-1$&0&$0$&1&0 \\\hline
7&0&$-1$&$0$&1&$-1$&1&8&$0$&1&0&$0$&0&$-1$\\\hline
9&0&$0$&$1$&$-1$&0&1&10&$0$&$-1$&0&$1$&0&$-1$\\\hline
11&1&$0$&$-1$&$0$&0&1&12&$0$&0&1&$-1$&1&$-1$\\\hline
13&0&0&1&0&$-1$&0&14&$-1$&0&0&0&0&1\\\hline
15&1&0&$-1$&0&1&$-1$&16&1&0&$-1$&1&$-1$&0\\\hline
17&$-1$&0&0&0&1&$-1$&18&1&1&$0$&$-1$&0&0\\\hline
19&$-1$&0&0&1&$-1$&0&20&1&$-1$&$0$&0&0&0\\\hline
21&$-1$&1&1&$-1$&0&0&22&0&1&$-1$&0&0&0\\\hline
23&$-1$&$-1$&1&0&0&0&24&0&$-1$&$-1$&1&0&0\\\hline
25&0&0&0&$-1$&1&0&26&0&0&0&$0$&$-1$&1\\\hline 27&0&0&0&0&0&$-1$
&&&&&&&\\\hline\end{tabular}\end{center}
$$E_{-\al_1}|_V=x_2\ptl_{x_1}-x_{14}\ptl_{x_{11}}-x_{17}\ptl_{x_{15}}
-x_{19}\ptl_{x_{16}}-x_{21}\ptl_{x_{18}}-x_{23}\ptl_{x_{20}},\eqno(5.21)$$
$$E_{-\al_2}|_V=x_6\ptl_{x_4}+x_7\ptl_{x_5}+x_{10}\ptl_{x_8}
-x_{20}\ptl_{x_{18}}-x_{23}\ptl_{x_{21}}-x_{24}\ptl_{x_{22}},\eqno(5.22)$$
$$E_{-\al_3}|_V=x_3\ptl_{x_2}-x_{11}\ptl_{x_9}
-x_{15}\ptl_{x_{12}}-x_{16}\ptl_{x_{13}}-x_{22}\ptl_{x_{21}}-x_{24}
\ptl_{x_{23}},\eqno(5.23)$$
$$E_{-\al_4}|_V=x_4\ptl_{x_3}+x_9\ptl_{x_7}+x_{12}\ptl_{x_{10}}
+x_{18}\ptl_{x_{16}}+x_{21}\ptl_{x_{19}}-x_{25}\ptl_{x_{24}},\eqno(5.24)$$
$$E_{-\al_5}|_V=x_5\ptl_{x_4}+x_7\ptl_{x_6}+x_{13}\ptl_{x_{12}}+
x_{16}\ptl_{x_{15}}
+x_{19}\ptl_{x_{17}}-x_{26}\ptl_{x_{25}},\eqno(5.25)$$
$$E_{-\al_6}|_V=x_8\ptl_{x_5}+x_{10}\ptl_{x_7}+x_{12}\ptl_{x_9}+x_{15}
\ptl_{x_{11}} +x_{17}\ptl_{x_{14}}-x_{27}\ptl_{x_{26}},\eqno(5.26)$$
(e.g., cf. [27]). Moreover,
$$E_{\al_r}(x_i)\neq 0\Leftrightarrow a_{r,i}<0,\;\;
E_{-\al_r}(x_i)\neq 0\Leftrightarrow a_{r,i}>0.\eqno(5.27)$$
\pse

{\bf Theorem 5.1}. {\it The ternary weight code ${\cal C}_{E_6,1}$
of $E_6$ on $V_{E_6}$ is an orthogonal $[27,6,12]$-code.}\psp

{\it Proof}. Write
$$A_{E_6}=(a_{r,i})_{6\times 27}.\eqno(5.28)$$
Denote by $\xi_r$ the $r$th row of the matrix $A_{E_6}$. Then
$$\wt\xi_r=12\qquad\for\;\;r\in\ol{1,6}.\eqno(5.29)$$
Moreover,
$$\wt(\xi_1+\xi_3)=\wt(\xi_2+\xi_4)=12,\;\;\wt (\xi_1+\xi_4)=18, \eqno(5.30)$$
$$\wt (\xi_1+\xi_2)=\wt(\xi_2+\xi_3)=\wt(\xi_2+\xi_5)=\wt(\xi_2+\xi_6)=18.\eqno(5.31)$$
By an argument as (3.29)-(3.31) and symmetry, we have
$$\xi_i\cdot\xi_j\equiv 0\;(\mbox{mod}\;3)\qquad\for\;\;i,j\in\ol{1,6},\eqno(5.32)$$
that is ${\cal C}_{E_6,1}$ is orthogonal.

 Note that
the Lie subalgebra ${\cal G}^{E_6}_{A,1}$ generated by
$\{E_{\pm\al_i}\mid 2\neq i\in\ol{1,6}\}$ is isomorphic to
$sl(6,\mbb{C})$. Recall that a singular vector in a module of
simple Lie algebra is a weight vector annihilated by its positive
root vectors. By Table 5.1 and (5.27), the ${\cal
G}^{E_6}_{A,1}$-singular vectors are $x_1$ of weight $\lmd_1$,
$x_6$ of  weight $\lmd_4$ and $x_{20}$ of weight $\lmd_1$. So the
$({\cal G}^{E_6},{\cal G}^{E_6}_{A,1})$-branch rule on $V_{E_6}$
is
$$V_{E_6}\cong V_{A_5}(\lmd_1)\oplus V_{A_5}(\lmd_4)\oplus V_{A_5}(\lmd_1).\eqno(5.33)$$
Denote by ${\cal G}^{E_6}_{A,2}$ the Lie subalgebra of ${\cal
G}^{E_6}$ generated by $\{E_{\pm\al_r},\;E_{\pm(\al_2+\al_4)}\mid
2,4\neq r\in\ol{1,6}\}$. The algebra ${\cal G}^{E_6}_{A,2}$ is also
isomorphic to $sl(6,\mbb{C})$. According to Table 5.1 and (5.27),
the ${\cal G}^{E_6}_{A,2}$-singular vectors are $x_1$ of  weight
$\lmd_1$, $x_4$ of  weight $\lmd_4$ and $x_{18}$ of  weight
$\lmd_1$. Hence (5.33) is also the $({\cal G}^{E_6},{\cal
G}^{E_6}_{A,2})$-branch rule. Since the module $V_{A_5}(\lmd_2)$ is
contragredient to $V_{A_5}(\lmd_4)$, they have the same ternary
weight code. By (2.39) and (2.43) with $n=6$, the minimal distances
of the subcodes $\sum_{2\neq i\in\ol{1,6}}\mbb{F}_3\xi_i$ and
$\mbb{F}_3(\xi_2+\xi_4)+\sum_{2,4\neq i\in\ol{1,6}}\mbb{F}_3\xi_i$
are $\wt\xi_1=12$.

Recall $\frac{1}{2}=-1$ in $\mbb{F}_3$. Moreover,
$$-(\al_2+\al_4)=-\ves_1-\ves_2+\ves_3-\ves_4+\ves_5+\ves_6+\sqrt{2}\ves_7\;\;\mbox{in}\;\;{\cal H}_{E_6,3}.\eqno(5.34)$$
Thus in ${\cal H}_{E_6,3}$,
$$\al_1-(\al_2+\al_4)=\ves_2+\ves_3-\ves_4+\ves_5+\ves_6+\sqrt{2}\ves_7,\eqno(5.35)$$
$$\al_1-\al_2-(\al_2+\al_4)=-\ves_3-\ves_4+\ves_5+\ves_6+\sqrt{2}\ves_7,\eqno(5.36)$$
$$\al_1-\al_2-(\al_2+\al_4)+\al_6=-\ves_3-\ves_4-\ves_5+\sqrt{2}\ves_7,\eqno(5.37)$$
$$\al_1-\al_2-(\al_2+\al_4)-\al_5+\al_6=-\ves_3+\ves_4+\sqrt{2}\ves_7.\eqno(5.38)$$
Note that
$$\wt (\xi_1-(\xi_2+\xi_4)),\;\wt(\xi_1-\xi_2-(\xi_2+\xi_4))\geq 12,\eqno(5.39)$$
$$\wt(\xi_1-\xi_2-(\xi_2+\xi_4)+\xi_6),\;\;\wt(\xi_1-\xi_2-(\xi_2+\xi_4)-\xi_5+\xi_6)\geq 12\eqno(5.40)$$
because the minimal distance of
$\mbb{F}_3(\xi_2+\xi_4)+\sum_{2,4\neq i\in\ol{1,6}}\mbb{F}_3\xi_i$
is 12. Furthermore,
$$-\sum_{i=1}^6\ves_i+\sqrt{2}\ves_7=\al_1-\al_2-\al_3\qquad\mbox{in}\;\;{\cal H}_{E_6,3}.\eqno(5.41)$$
We calculate
$$\wt(\xi_1-\xi_2-\xi_3)=21.\eqno(5.42)$$
By (5.13), the minimal distance of the ternary code ${\cal
C}_{E_6,1}$ is 12. $\qquad\Box$\psp

Next we consider the ternary weight code ${\cal C}_{E_6,2}$ of $E_6$
on its adjoint module. Take any order
$$\{y_1,...,y_{36}\}=\{E_\al\mid \al\in \Phi^+_{E_6}\}.\eqno(5.43)$$
Write
$$[\al_i,y_j]=b_{i,j},\qquad B_{E_6}=(b_{i,j})_{6\times 36}.\eqno(5.44)$$
\pse

{\bf Theorem 5.2}. {\it The ternary weight code ${\cal C}_{E_6,2}$
(generated $B_{E_6}$) of $E_6$ on its adjoint module is an
orthogonal $[36,5,21]$-code.}

{\it Proof}. Denote by $\zeta_i$ the $i$th row of $B_{E_6}$. Note
that
$$\zeta_1-\zeta_3+\zeta_5-\zeta_6\equiv 0\qquad\mbox{in}\;\;\mbb{F}_3.\eqno(5.45)$$
Thus
$${\cal C}_{E_6,2}=\sum_{i=2}^6\mbb{F}_3\zeta_i.\eqno(5.46)$$
Denote by ${\cal G}^{E_6}_D$
the Lie subalgebra  of ${\cal G}^{E_6}$  generated by $\{E_{\pm \al_r}\mid
r\in\ol{2,6}\}$. According to the Dynkin diagram of $E_6$,
$${\cal G}^{E_6}_D\cong o(10,\mbb{C}).\eqno(5.47)$$
Let ${\cal G}^{E_6}_+=\sum_{i=1}^{36}\mbb{C}y_i$ and denote by
${\cal G}^{E_6}_{D,+}$ the subspace spanned by the  root vectors
$E_\al\in{\cal G}^{E_6}_D$ with $\al\in\Phi_{E_6}^+$. Then $[{\cal
G}^{E_6}_{D,+},{\cal G}^{E_6}_+]\subset {\cal G}^{E_6}_+$. Moreover,
 the space
${\cal G}^{E_6}_+$ contains  ${\cal G}^{E_6}_D$-singular vectors
$E_{\al_4+\al_5+\sum_{i=2}^6\al_i}$ of weight $\lmd_2$ (the highest
root) and $E_{\al_2+\al_4+\sum_{r=3}^5\al_r+\sum_{i=1}^6\al_i}$ of
weight $\lmd_5$. Hence, we have the partial $({\cal G}_{E_6},{\cal
G}^{E_6}_D)$-branch rule on ${\cal G}_{E_6}$:
$${\cal G}_{E_6}^+\cong{\cal G}^{E_6}_{D+}\oplus V_{D_5}(\lmd_5).
\eqno(5.48)$$ Thus the ternary weight code ${\cal C}_{E_6,2}$ of
$E_6$ on its adjoint module  is exactly the code ${\cal C}_2\oplus
{\cal C}_3({\cal V})$ with $n=5$ in Corollary 3.4, which is a
ternary orthogonal $[36,5,21]$-code.$\qquad\Box$

\section{Representations of $E_7, E_8$ and Ternary Codes}

In this section, we study the ternary weight codes of $E_7$ on its
minimal irreducible module and adjoint module, and the ternary
weight code of $E_8$ on its minimal irreducible module (adjoint
module).

Recall the settings in (2.2)-(2.4) and (5.1)-(5.8). Taking $n=8$, we
have the root system of $E_7$:
 $$\Phi_{E_7}=\left\{\ves_i-\ves_j,\frac{1}{2}\sum_{s=1}^8\iota_s\ves_s
\mid i,j\in\ol{1,8},\;i\neq j;\;\iota_s=\pm
1,\;\sum_{s=1}^8\iota_s=0\right\}\eqno(6.1)$$ and the simple
positive roots are:
$$\al_1=\ves_2-\ves_3,\;\al_2=\frac{1}{2}\sum_{j=1}^4(\ves_{4+j}-\ves_j),\;\;\al_i=\ves_i-\ves_{i+1},\qquad i\in\ol{3,7}.
\eqno(6.2)$$

The Dynkin diagram of $E_7$ is as follows:

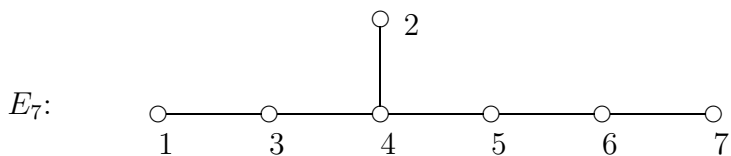
\begin{picture}(93,20)
\put(2,0){$E_7$:}\put(21,0){\circle{2}}\put(21,
-5){1}\put(22,0){\line(1,0){12}}\put(35,0){\circle{2}}\put(35,
-5){3}\put(36,0){\line(1,0){12}}\put(49,0){\circle{2}}\put(49,
-5){4}\put(49,1){\line(0,1){10}}\put(49,12){\circle{2}}\put(52,10){2}\put(50,0){\line(1,0){12}}
\put(63,0){\circle{2}}\put(63,-5){5}\put(64,0){\line(1,0){12}}\put(77,0){\circle{2}}\put(77,
-5){6}\put(78,0){\line(1,0){12}}\put(91,0){\circle{2}}\put(91,
-5){7}
\end{picture}
\vspace{0.7cm}

The minimal module $V_{E_7}$ of $E_7$ is of  56-dimensional and
has a basis $\{x_i\mid i\in\ol{1,56}\}$ with the representation
formulas determined by
\begin{eqnarray*}E_{\al_1}|_V&=&-x_6\ptl_{x_8}-x_9\ptl_{x_{11}}
-x_{10}\ptl_{x_{13}}-x_{12}\ptl_{x_{16}}-x_{14}\ptl_{x_{19}}
-x_{17}\ptl_{x_{22}}
\\ &&+x_{35}\ptl_{x_{40}}+x_{38}\ptl_{x_{43}}+x_{41}\ptl_{x_{45}}
+x_{44}\ptl_{x_{47}}+x_{46}\ptl_{x_{48}}+x_{49}\ptl_{x_{51}},\hspace{2.6cm}(6.3)
\end{eqnarray*}
\begin{eqnarray*}E_{\al_2}|_V&=&x_5\ptl_{x_7}+x_6\ptl_{x_9}
+x_8\ptl_{x_{11}}-x_{20}\ptl_{x_{23}}-x_{24}\ptl_{x_{26}}-x_{27}\ptl_{x_{29}}
\\ &&-x_{28}\ptl_{x_{30}}-x_{31}\ptl_{x_{33}}-x_{34}\ptl_{x_{37}}
+x_{46}\ptl_{x_{49}}+x_{48}\ptl_{x_{51}}+x_{50}\ptl_{x_{52}},\hspace{2.6cm}(6.4)
\end{eqnarray*}
\begin{eqnarray*}E_{\al_3}|_V&=&-x_5\ptl_{x_6}-x_7\ptl_{x_9}
-x_{13}\ptl_{x_{15}}-x_{16}\ptl_{x_{18}}-x_{19}\ptl_{x_{21}}
-x_{22}\ptl_{x_{25}}
\\ &&+x_{32}\ptl_{x_{35}}+x_{36}\ptl_{x_{38}}+x_{39}\ptl_{x_{41}}
+x_{42}\ptl_{x_{44}}+x_{48}\ptl_{x_{50}}+x_{51}\ptl_{x_{52}},
\hspace{2.6cm}(6.5)
\end{eqnarray*}
\begin{eqnarray*}E_{\al_4}|_V&=&x_4\ptl_{x_5}-x_9\ptl_{x_{10}}
-x_{11}\ptl_{x_{13}}-x_{18}\ptl_{x_{20}}-x_{21}\ptl_{x_{24}}
-x_{25}\ptl_{x_{28}}
\\ &&-x_{29}\ptl_{x_{32}}-x_{33}\ptl_{x_{36}}-x_{37}\ptl_{x_{39}}
-x_{44}\ptl_{x_{46}}-x_{47}\ptl_{x_{48}}+x_{52}\ptl_{x_{53}},
\hspace{2.6cm}(6.6)
\end{eqnarray*}
\begin{eqnarray*}E_{\al_5}|_V&=&x_3\ptl_{x_4}-x_{10}\ptl_{x_{12}}
-x_{13}\ptl_{x_{16}}-x_{15}\ptl_{x_{18}}-x_{24}\ptl_{x_{27}}
-x_{26}\ptl_{x_{29}}
\\ &&-x_{28}\ptl_{x_{31}}-x_{30}\ptl_{x_{33}}-x_{39}\ptl_{x_{42}}
-x_{41}\ptl_{x_{44}}-x_{45}\ptl_{x_{47}}+x_{53}\ptl_{x_{54}},
\hspace{2.6cm}(6.7)
\end{eqnarray*}
\begin{eqnarray*}E_{\al_6}|_V&=&x_2\ptl_{x_3}-x_{12}\ptl_{x_{14}}
-x_{16}\ptl_{x_{19}}-x_{18}\ptl_{x_{21}}-x_{20}\ptl_{x_{24}}
-x_{23}\ptl_{x_{26}}
\\ &&-x_{31}\ptl_{x_{34}}-x_{33}\ptl_{x_{37}}-x_{36}\ptl_{x_{39}}
-x_{38}\ptl_{x_{41}}-x_{43}\ptl_{x_{45}}+x_{54}\ptl_{x_{55}},
\hspace{2.6cm}(6.8)
\end{eqnarray*}
\begin{eqnarray*}E_{\al_7}|_V&=&x_1\ptl_{x_2}-x_{14}\ptl_{x_{17}}
-x_{19}\ptl_{x_{22}}-x_{21}\ptl_{x_{25}}-x_{24}\ptl_{x_{28}}-x_{26}\ptl_{x_{30}}
\\ &&-x_{27}\ptl_{x_{31}}-x_{29}\ptl_{x_{33}}-x_{32}\ptl_{x_{36}}
-x_{35}\ptl_{x_{38}}-x_{40}\ptl_{x_{43}}+x_{55}\ptl_{x_{56}},
\hspace{2.6cm}(6.9)
\end{eqnarray*}
\begin{eqnarray*}E_{-\al_1}|_V&=&x_8\ptl_{x_6}+x_{11}\ptl_{x_9}
+x_{13}\ptl_{x_{10}}+x_{16}\ptl_{x_{12}}+x_{19}\ptl_{x_{14}}+x_{22}\ptl_{x_{17}}
\\&&-x_{40}\ptl_{x_{35}}-x_{43}\ptl_{x_{38}}-x_{45}\ptl_{x_{41}}
-x_{47}\ptl_{x_{44}}-x_{48}\ptl_{x_{46}}-x_{51}\ptl_{x_{49}},
\hspace{2.2cm}(6.10)
\end{eqnarray*}
\begin{eqnarray*}E_{-\al_2}|_V&=&-x_7\ptl_{x_5}-x_9\ptl_{x_6}
-x_{11}\ptl_{x_8}+x_{23}\ptl_{x_{20}}+x_{26}\ptl_{x_{24}}+
x_{29}\ptl_{x_{27}}\\
&& +x_{30}\ptl_{x_{28}}+x_{33}\ptl_{x_{31}}+x_{37}\ptl_{x_{34}}
-x_{49}\ptl_{x_{46}}-x_{51}\ptl_{x_{48}}-x_{52}\ptl_{x_{50}},
\hspace{2.2cm}(6.11)
\end{eqnarray*}
\begin{eqnarray*}E_{-\al_3}|_V&=&x_6\ptl_{x_5}+x_9\ptl_{x_7}
+x_{15}\ptl_{x_{13}}+x_{18}\ptl_{x_{16}}+x_{21}\ptl_{x_{19}}+x_{25}\ptl_{x_{22}}
\\&&-x_{35}\ptl_{x_{32}}-x_{38}\ptl_{x_{36}}-x_{41}\ptl_{x_{39}}
-x_{44}\ptl_{x_{42}}-x_{50}\ptl_{x_{48}}-x_{52}\ptl_{x_{51}},
\hspace{2.2cm}(6.12)
\end{eqnarray*}
\begin{eqnarray*}E_{-\al_4}|_V&=&-x_5\ptl_{x_4}+x_{10}\ptl_{x_9}
+x_{13}\ptl_{x_{11}}+x_{20}\ptl_{x_{18}}+x_{24}\ptl_{x_{21}}+x_{28}\ptl_{x_{25}}
\\&&+x_{32}\ptl_{x_{29}}+x_{36}\ptl_{x_{33}}+x_{39}\ptl_{x_{37}}
+x_{46}\ptl_{x_{44}}+x_{48}\ptl_{x_{47}}-x_{53}\ptl_{x_{52}},
\hspace{2.2cm}(6.13)
\end{eqnarray*}
\begin{eqnarray*}E_{-\al_5}|_V&=&-x_4\ptl_{x_3}+x_{12}\ptl_{x_{10}}
+x_{16}\ptl_{x_{13}}+x_{18}\ptl_{x_{15}}+x_{27}\ptl_{x_{24}}
+x_{29}\ptl_{x_{26}}\\&&
+x_{31}\ptl_{x_{28}}+x_{33}\ptl_{x_{30}}+x_{42}\ptl_{x_{39}}
+x_{44}\ptl_{x_{41}}+x_{47}\ptl_{x_{45}}-x_{54}\ptl_{x_{53}},
\hspace{2.2cm}(6.14)
\end{eqnarray*}
\begin{eqnarray*}E_{-\al_6}|_V&=&-x_3\ptl_{x_2}+x_{14}\ptl_{x_{12}}
+x_{19}\ptl_{x_{16}}+x_{21}\ptl_{x_{18}}+x_{24}\ptl_{x_{20}}+x_{26}\ptl_{x_{23}}
\\&&+x_{34}\ptl_{x_{31}}+x_{37}\ptl_{x_{33}}+x_{39}\ptl_{x_{36}}
+x_{41}\ptl_{x_{38}}+x_{45}\ptl_{x_{43}}-x_{55}\ptl_{x_{54}},
\hspace{2.2cm}(6.15)
\end{eqnarray*}
\begin{eqnarray*}E_{-\al_7}|_V&=&-x_2\ptl_{x_1}+x_{17}\ptl_{x_{14}}
+x_{22}\ptl_{x_{19}}+x_{25}\ptl_{x_{21}}+x_{28}\ptl_{x_{24}}
+x_{30}\ptl_{x_{26}}\\&&
+x_{31}\ptl_{x_{27}}+x_{33}\ptl_{x_{29}}+x_{36}\ptl_{x_{32}}
+x_{38}\ptl_{x_{35}}+x_{43}\ptl_{x_{40}}-x_{56}\ptl_{x_{55}},
\hspace{2.2cm}(6.16)
\end{eqnarray*}
$$h_r|_V=\sum_{i=1}^{28}a_{r,i}(x_i\ptl{x_i}-x_{57-i}\ptl{x_{57-
i}})\qquad\for\;\;r\in\ol{1,7},\eqno(6.17)$$ where $a_{r,i}$ are
constants given by the following table:
\begin{center}{\bf \large Table 6.1}\end{center}
\begin{center}\begin{tabular}{|r||r|r|r|r|r|r|r||r||r|r|r|r|r|r|r|}\hline
$i$&$a_{1,i}$&$a_{2,i}$&$a_{3,i}$&$a_{4,i}$&$a_{5,i}$&$a_{6,i}$&$a_{7,i}$
&$i$&$a_{1,i}$&$a_{2,i}$&$a_{3,i}$&$a_{4,i}$&$a_{5,i}$&
$a_{6,i}$&$a_{7,i}$
\\\hline\hline 1&0&0&0&0&0&0&1&2&0&0&0&0&0&1&$-1$\\\hline
3&0&0&0&0&1&$-1$&0&4&$0$&0&0&1&$-1$&0&0
\\\hline 5&0&1&1&$-1$&0&0&0&6&1&1&$-1$&0&0&0&0 \\\hline
7&0&$-1$&$1$&0&0&0&0&8&$-1$&1&0&$0$&0&0&0\\\hline
9&1&$-1$&$-1$&1&0&0&0&10&1&0&0&$-1$&1&0&0\\\hline
11&$-1$&$-1$&0&1&$0$&0&0&12&1&0&0&0&$-1$&1&0\\\hline
13&$-1$&0&1&$-1$&1&$0$&0&14&1&0&0&0&0&$-1$&1\\\hline
15&0&0&$-1$&0&1&0&0&16&$-1$&0&1&0&$-1$&1&0\\\hline
17&1&0&0&0&0&0&$-1$&18&0&0&$-1$&1&$-1$&1&0\\\hline
19&$-1$&0&1&0&0&$-1$&1&20&0&1&0&$-1$&0&1&0\\\hline
21&0&0&$-1$&1&0&$-1$&1&22&$-1$&0&1&0&0&0&$-1$\\\hline
23&0&$-1$&0&0&0&1&0&24&0&1&0&$-1$&1&$-1$&1\\\hline
25&0&0&$-1$&1&0&0&$-1$&26&0&$-1$&0&0&1&$-1$&1\\\hline
27&0&1&0&0&$-1$&0&1&28&0&1
&0&$-1$&1&0&$-1$\\\hline\end{tabular}\end{center} (e.g., cf. [28]).
 Again we have
$$E_{\al_r}(x_i)\neq 0\Leftrightarrow a_{r,i}<0,\;\;
E_{-\al_r}(x_i)\neq 0\Leftrightarrow a_{r,i}>0.\eqno(6.18)$$
Denote
$$A_{E_7}=(a_{r,i})_{7\times 28}.\eqno(6.19)$$
\pse

{\bf Theorem 6.1}. {\it The ternary weight code ${\cal C}_{E_7,1}$ of $E_7$ on $V_{E_7}$ is an
orthogonal $[28,7,12]$-code}.

{\it Proof}.
Note that the root system of $A_7$:
$$\Phi_{A_7}=\{\ves_i-\ves_j\mid i,j\in\ol{1,8},\;i\neq j\}\subset \Phi_{E_7}.\eqno(6.20)$$
Thus we have the Lie subalgebra of ${\cal G}^{E_7}$ (cf. (5.1)-(5.7)
with $m=7$):
$${\cal G}^{E_7}_A=\sum_{i=1}^7\mbb{C}\al_i+\sum_{\al\in\Phi_{A_7}}\mbb{C}E_\al\cong sl(8,\mbb{C})
.\eqno(6.21)$$
Moreover,
$$\al_1'=\ves_1-\ves_2=-2\al_2-2\al_1-3\al_3-4\al_4-3\al_5-2\al_6-\al_7.\eqno(6.22)$$
Note that $x_{23}$ is a ${\cal G}^{E_7}_A$-singular vector of weight
$\lmd_6$ and $x_{49}$ is a ${\cal G}^{E_7}_A$-singular vector of
weight $\lmd_2$ by (6.17), (6.18) and Table 6.1. Thus the $({\cal
G}^{E_7},{\cal G}^{E_7}_A)$-branch rule on $V_{E_7}$ is
$$V_{E_7}\cong V_{A_7}(\lmd_2)\oplus V_{A_7}(\lmd_6).\eqno(6.23)$$
Since  $V_{A_7}(\lmd_6)$ is contragredient to $V_{A_7}(\lmd_2)$,
they have the same ternary weight code of ${\cal G}^{E_7}_A$, which
is the ${\cal C}_3({\cal A}_2)$ with $m=2$ in Theorem 2.3. Hence the
weight matrix of ${\cal G}^{E_7}_A$ on $V_{E_7}$ generates a ternary
orthogonal $[56,7,24]$-code.

On the other hand,
$$\sum_{i=1}^7\mbb{F}_3\al_i=\mbb{F}_3\al_1'+\sum_{2\neq i\in\ol{1,7}}\mbb{F}_3\al_i\eqno(6.24)$$
by (6.1) and  the fact $1/2\equiv -1$ in $\mbb{F}_3$. Thus the
weight matrix $(A_{E_7},-A_{E_7})$ of $E_7$ on $V_{E_7}$ generates
the same ternary code as the weight matrix of ${\cal G}^{E_7}_A$
on $V_{E_7}$. So $(A_{E_7},-A_{E_7})$ generates a ternary
orthogonal $[56,7,24]$-code. Hence the ternary code ${\cal
C}_{E_7,1}$ generated by $A_{E_7}$ is  an orthogonal
$[28,7,12]$-code.$\qquad\Box$\psp

Next we consider the ternary weight code of $E_7$ on its adjoint
module. Recall the construction of ${\cal G}^{E_7}$ in (5.1)-(5.7)
with $m=7$. The  $({\cal G}^{E_7},{\cal G}^{E_7}_A)$-branch rule
on ${\cal G}^{E_7}$ is
$${\cal G}^{E_7}\cong {\cal G}^{E_7}_A\oplus V_{A_7}(\lmd_4).\eqno(6.25)$$
The module $V_{A_7}(\lmd_4)$ of $sl(8,\mbb{C})\;(\cong {\cal
G}^{E_7}_A)$ is exactly ${\cal A}_4$ in (2.10) with $n=8$, which is
self-dual. For convenience, we study the ternary code generated by
the weight matrix of $sl(8,\mbb{C})$ on ${\cal A}_4$. Taking any
order of its basis
$$\{z_1,...,z_{70}\}=\{\sta_{i_1}\sta_{i_2}\sta_{i_3}\sta_{i_4}\mid 1\leq i_1<i_2<i_3<i_4\leq 8\},\eqno(6.26)$$
we write
$$[E_{r,r},z_i]=b_{r,i}z_i,\qquad B_{E_7}=(b_{r,i})_{7\times 70}.\eqno(6.27)$$
Denote by $\eta_r$ the $r$th row of $B_{E_7}$ and by ${\cal C}'$ the ternary code generated by $B_{E_7}$.
Set
$$v(s,t)=\sum_{i=1}^s\eta_i-\sum_{j=1}^t\eta_{s+j}\qquad\in {\cal C}'.\eqno(6.28)$$
Moreover, we only calculate the related weights:
\begin{center}{\bf Table 6.2}\end{center}
\begin{center}\begin{tabular}{|c|c|c|c|c|c|c|c|c|}\hline
(s,t)&(1,1)&(2,2)&(3,3)&(4,4)&(3,0)&(6,0)&(4,1)&(5,2)
\\\hline \wt v(s,t)&40&44&48&34&60&30&46&50\\\hline\end{tabular}\end{center}
Recall (2.65)-(2.70). We have
\begin{center}{\bf Table 6.3}\end{center}
\begin{center}\begin{tabular}{|c|c|c|c|c|c|c|c|c|}\hline
(s,t)&(1,1)&(2,2)&(3,3)&(4,4)&(3,0)&(6,0)&(4,1)&(5,2)
\\\hline 2\wt u(s,t)&26&40&42&32&30&24&38&34\\\hline\end{tabular}\end{center}
According to (6.1), the Weyl group ${\cal W}_{E_7}$ contains the
permutation group $S_8$ on the sub-indices of $\ves_i$. By (1.9),
(1.11) and the values of $\wt v(s,t)+2\wt u(s,t)$ from the above
tables, 54, 66, 84 and 90 are the only weights of the nonzero
codewords in ${\cal C}_3({\cal G}^{E_7})$, the ternary code
generated by the weight matrix of ${\cal G}^{E_7}_A$ on ${\cal
G}^{E_7}$. By (6.24) and an argument as (3.29)-(3.31), we have:
\psp

{\bf Theorem 6.2}. {\it The ternary weight code of $E_7$ on its adjoint module is an orthogonal
$[63,7,27]$-code}.\psp

The minimal representation of $E_8$ is its adjoint module. Recall
the settings in (2.2)-(2.4) and
 construction of the simple Lie algebra ${\cal G}^{E_8}$ given in (5.1)-(5.8) with $m=8$. we have the $E_8$ root system
$$\Phi_{E_8}=\left\{\pm \ves_i\pm \ves_j,\frac{1}{2}\sum_{i=1}^8\iota_i\ves_i
\mid i,j\in\ol{1,8},\;i\neq j;\;\iota_i=\pm
1,\;\sum_{i=1}^8\iota_i\in 2\mbb{Z}\right\}\eqno(6.29)$$
and positive simple roots:
$$\al_1=\frac{1}{2}(\sum_{j=2}^7\ves_j-\ves_1-\ves_8),\;\al_2=-\ves_1-\ves_2,\;\al_r=\ves_{r-2}-\ves_{r-1},\qquad
r\in\ol{3,8}.\eqno(6.30)$$
The Dynkin diagram of $E_8$ is as
follows:

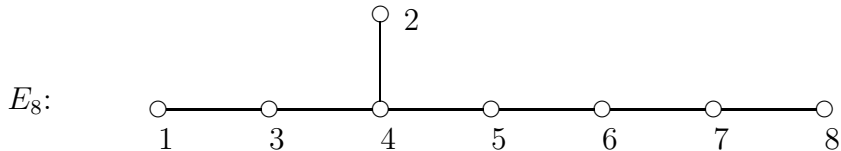
\begin{picture}(110,23)
\put(2,0){$E_8$:}\put(21,0){\circle{2}}\put(21,
-5){1}\put(22,0){\line(1,0){12}}\put(35,0){\circle{2}}\put(35,
-5){3}\put(36,0){\line(1,0){12}}\put(49,0){\circle{2}}\put(49,
-5){4}\put(49,1){\line(0,1){10}}\put(49,12){\circle{2}}\put(52,10){2}\put(50,0){\line(1,0){12}}
\put(63,0){\circle{2}}\put(63,-5){5}\put(64,0){\line(1,0){12}}\put(77,0){\circle{2}}\put(77,
-5){6}\put(78,0){\line(1,0){12}}\put(91,0){\circle{2}}\put(91,
-5){7}\put(92,0){\line(1,0){12}}\put(105,0){\circle{2}}\put(105,
-5){8}
\end{picture}
\vspace{0.7cm}

Observe that the root system of $o(16,\mbb{C})$:
$$\Phi_{D_8}=\{\pm \ves_i\pm \ves_j
\mid i,j\in\ol{1,8},\;i\neq j\}\subset \Phi_{E_8}.\eqno(6.31)$$
So the Lie subalgebra
$${\cal G}^{E_8}_D=H_{E_8}+\sum_{\al\in\Phi_{D_8}}\mbb{C}E_\al\eqno(6.32)$$
of ${\cal G}^{E_8}$ is exactly isomorphic to $o(16,\mbb{C})$.
Moreover, the $({\cal G}^{E_8},{\cal G}^{E_8}_D)$-branch rule on
${\cal G}^{E_8}$ is
$${\cal G}^{E_8}\cong {\cal G}^{E_8}_D\oplus V_{D_8}(\lmd_8).\eqno(6.33)$$
In fact, $V_{D_8}(\lmd_8)$ is exactly the spin module ${\cal V}$ in
(3.35). Since
$$\sum_{i=1}^8\mbb{F}_3\al_i=\sum_{\al\in\Phi_{D_8}}\mbb{F}_3\al,\eqno(6.34)$$
the ternary weight code of $E_8$ on ${\cal G}^{E_8}$ is the same as that of ${\cal G}^{E_8}_D$
on ${\cal G}^{E_8}$. By Corollary 3.4 with $m=8$, we have:
\psp

{\bf Theorem 6.3}. {\it The ternary weight code of $E_8$ on its adjoint module is an orthogonal $[120,8,57]$-code
}.
\vspace{1cm}
\begin{center}{\Large \bf Acknowledgments}\end{center}

Part of this work was done when the author visited The University of Sydney. He would like to thank Prof. Ruibin Zhang
for his hospitality and the enthusiastic academic discussions.

\bibliographystyle{amsplain}

\end{document}